\documentclass[12pt]{amsart}
\usepackage[T1]{fontenc}
\usepackage{lmodern}
\usepackage[margin=3cm]{geometry}
\usepackage{amsmath,amssymb,amsthm}
\usepackage{booktabs}
\usepackage{needspace}
\usepackage{microtype}
\usepackage{xcolor}
\definecolor{linkblue}{rgb}{0,0.2,0.55}
\usepackage[colorlinks=true,linkcolor=linkblue,citecolor=linkblue,urlcolor=linkblue,
  pdftitle={On mod p Hecke algebras of abstract p-adic automorphic forms},
  pdfauthor={Nha Xuan Truong},
  pdfsubject={Number theory; mod p Hecke algebras; p-adic automorphic forms},
  pdfkeywords={mod p Hecke algebra, nilpotence method, p-power level, completed homology}]{hyperref}

\theoremstyle{plain}
\newtheorem{theorem}{Theorem}[section]
\newtheorem{lemma}[theorem]{Lemma}
\newtheorem{proposition}[theorem]{Proposition}
\newtheorem{corollary}[theorem]{Corollary}
\theoremstyle{definition}
\newtheorem{definition}[theorem]{Definition}
\newtheorem{remark}[theorem]{Remark}
\newtheorem{example}[theorem]{Example}

\theoremstyle{plain}

\theoremstyle{plain}
\newtheorem*{theorem*}{Theorem}
\newtheorem*{lemma*}{Lemma}
\newtheorem*{proposition*}{Proposition}
\newtheorem*{corollary*}{Corollary}
\newtheorem*{corABC*}{Corollary}
\theoremstyle{definition}
\newtheorem*{definition*}{Definition}
\newtheorem*{remark*}{Remark}
\newtheorem*{hypothesis*}{Hypothesis}
\numberwithin{equation}{section}

\newcommand{\F}{\mathbb{F}}
\newcommand{\Q}{\mathbb{Q}}
\newcommand{\Z}{\mathbb{Z}}
\newcommand{\A}{\mathbb{A}}
\newcommand{\T}{\mathbb{T}}
\newcommand{\cO}{\mathcal{O}}
\newcommand{\Iw}{\mathrm{Iw}_p}
\newcommand{\Iwone}{\mathrm{Iw}^1_p}

\newcommand{\GL}{\mathrm{GL}}

\newcommand{\Sbar}{\bar S}
\newcommand{\Ht}{\widetilde H}
\newcommand{\m}{\mathfrak m}
\newcommand{\n}{\mathfrak n}

\newcommand{\ad}{\mathrm{ad}}
\newcommand{\Frob}{\mathrm{Frob}}

\DeclareMathOperator{\Hom}{Hom}
\DeclareMathOperator{\End}{End}
\DeclareMathOperator{\Ann}{Ann}
\DeclareMathOperator{\Span}{Span}
\DeclareMathOperator{\gr}{gr}
\DeclareMathOperator{\Map}{Map}
\DeclareMathOperator{\Ind}{Ind}
\DeclareMathOperator{\im}{im}

\newcommand{\smat}[4]{\left(\begin{smallmatrix}#1&#2\\#3&#4\end{smallmatrix}\right)}
\newcommand{\pmat}[4]{\begin{pmatrix}#1&#2\\#3&#4\end{pmatrix}}

\title[Mod $p$ Hecke algebras of $p$-power level]{On mod $p$ Hecke algebras of abstract $p$-adic automorphic forms}
\author{Nha Xuan Truong}
\address{Beijing International Center for Mathematical Research, Peking University,
5 Yi He Yuan Road, Haidian District, Beijing 100871, China}
\email{nxtruong@bicmr.pku.edu.cn}
\begin{document}
\raggedbottom

\begin{abstract}
We study completed local shallow Hecke algebras acting on residual
automorphic forms attached to projective augmented modules for
$\GL_2(\Q_p)$.
For odd primes $p$, we prove that the maximal ideal has nilpotence exponent
$O(p^{r/2})$ on the level-$r$ generalized eigenspace.
For $p\ge5$, if the corresponding Hecke eigenspace in the residual
smooth dual has finite length as a $\GL_2(\Q_p)$-representation,
then the completed local algebra is noetherian of Krull dimension
at least two.
If, in addition, its cotangent space has dimension at most two,
then it is a formal power series ring in two variables.
The nilpotence bound follows by conjugating level subgroups.
To deduce the dimension bound, we replace the Fourier pairing
in the classical nilpotence method with an inequality comparing
the dimension of each finite-level generalized eigenspace
with that of its Hecke eigenspace.
We apply these results to modular curves and definite quaternion algebras.
\end{abstract}

\maketitle

\section{Introduction}
\subsection{Historical background}\label{sec:intro-history}

Let $p$ be a prime and $\F$ a finite field of characteristic $p$.
Jochnowitz \cite{Jo82} studied the Hecke algebra generated by $T_\ell$,
$\ell\ne p$, on modular forms modulo $p$ of level one and all weights.
It decomposes into local factors indexed by residual Hecke eigensystems.
For $p=2$ it is local, and Nicolas and Serre \cite{NS12a,NS12b} proved
that it is a power series ring in $T_3$ and $T_5$.
For $p\ge5$, Bella\"iche and Khare \cite{BK15} proved that each local
factor has Krull dimension at least two and is a power series ring in
two variables when the residual representation is unobstructed.
Deo \cite{De17} extended these results to tame level $\Gamma_0(N)$ and
$\Gamma_1(N)$, with additional hypotheses for the power series
presentation at $\Gamma_1(N)$-level.
These dimension bounds use the characteristic-zero infinite fern of
Gouv\^ea and Mazur \cite{GM98}: a suitable Hecke algebra of $p$-adic
modular forms has dimension at least four, and specialization modulo
$p$ and a weight parameter gives the mod $p$ lower bound.

Medvedovsky \cite{Me15,Me18} developed an idea of Bella\"iche
\cite[Appendix]{BK15} into the \emph{nilpotence method}, which proves
dimension bounds entirely in characteristic $p$ and completes the
case $p=3$ begun in \cite{BK15}.
For a local factor $(A,\n)$, sublinear bounds in the weight for Hecke
nilpotence produce many forms killed by $\n^q$.
The perfect pairing $\langle T,f\rangle=a_1(Tf)$ between $A$ and
$\ker U_p$ converts
this count into a lower bound for
$\dim_\F A/\n^q$, and hence for $\dim A$.

The required estimates come from a nilpotence growth theorem for
degree-lowering recursion operators on $\F[y]$. This polynomial model is available when $X_0(p)$ has
genus zero, namely for $p\in\{2,3,5,7,13\}$.
Extending the argument to other primes would require a version of the growth
theorem for more general algebras \cite{Me15,Me18}.

We apply the nilpotence method to residual automorphic forms of
increasing $p$-power level.
The $\GL_2(\Q_p)$-action gives the nilpotence estimate, while local
representation theory bounds the eigenspaces.
This approach works directly with the level filtration and does not
require a polynomial presentation of classical mod $p$ forms.

\subsection{The Hecke algebra of \texorpdfstring{$p$}{p}-power level}
\label{sec:setting}\label{sec:algebra}

We work with the abstract $p$-adic automorphic forms of
Liu--Truong--Xiao--Zhao \cite[\S2]{LTXZ26}.
Let $\cO$ be the ring of integers of a finite extension of $\Q_p$,
with uniformizer $\varpi$ and residue field $\F$, and put
$G=\GL_2(\Q_p)$.
Write $\Iw$ for the standard Iwahori subgroup of $\GL_2(\Z_p)$,
consisting of matrices that are upper triangular modulo $p$.
Let $T(\F_p)\simeq\F_p^\times\times\F_p^\times$ be the diagonal
torus over $\F_p$.
Our input is a \emph{projective augmented module} $\Ht$: a compact
right $\cO[[\Iw]]$-module, finitely generated and projective over the
Iwahori group ring, with a compatible continuous right $G$-action.
The arithmetic examples are completed homology of modular curves,
localized at a non-Eisenstein maximal ideal, and completed homology
of a definite quaternion algebra split at $p$
(Examples \ref{ex:quaternion} and \ref{ex:modular} ).

Fix a character $\varepsilon:T(\F_p)\to\F^\times$.
For integral coefficients, we use its Teichm\"uller lift
$\widetilde\varepsilon:T(\F_p)\to\cO^\times$.
Let $w$ be the formal weight variable.
The associated space of abstract $p$-adic automorphic forms is
\[
S^{(\varepsilon)}_{p\text{-adic}}(\Ht)
=\Hom_{\cO[[\Iw]]}
\bigl(\Ht,C^0(\Z_p,\cO[[w]]^{(\varepsilon)})\bigr),
\]
with the twisted coefficient action of \cite[\S2.4]{LTXZ26}.
We study its specialization at $(\varpi,w)$:
\[
\Sbar=\Sbar^{(\varepsilon)}(\Ht)
=\Hom_{\F[[\Iw]]}\bigl(\Ht/\varpi,C^0(\Z_p,\F)\bigr).
\]
For simplicity, we suppress the character twist $\varepsilon$ on the
coefficient space. Explicitly,
\[
(h\|_u)(z)=\varepsilon(\bar a,\bar d)
h\Bigl(\frac{az+b}{cz+d}\Bigr),\qquad u=\smat abcd\in\Iw.
\]
Let $V_r\subseteq\Sbar$ consist of maps whose values factor through
$\Z_p/p^r\Z_p$. These spaces form the \emph{level filtration}
\[
V_0\subseteq V_1\subseteq\cdots,\qquad \Sbar=\bigcup_{r\ge0}V_r.
\]

Let $\T$ be a commutative $\F$-algebra of continuous endomorphisms
of $\Ht/\varpi$ commuting with $G$.
The action of $\T$ on $\Sbar$ is given by
\[
(t\phi)(x)=\phi(t(x))
\text{ for } t\in\T, \ \phi\in\Sbar,\text{ and } x\in\Ht/\varpi.\]
This action preserves each $V_r$.
In the arithmetic examples, $\T$ is generated by the unramified
Hecke operators $T_\ell$ and the central Hecke operators $S_\ell$
for primes $\ell$ outside a finite set $\mathcal R$ containing $p$.
We call this algebra \emph{shallow}: operators at the excluded
primes, including $U_p$, are omitted.
The set $\mathcal R$ is specified in each arithmetic example.

A maximal ideal $\m\subseteq\T$ \emph{occurs} if the localization
$\Sbar_\m$ is nonzero.
Let $V_{r,\m}$ be the generalized Hecke eigenspace in $V_r$
associated with $\m$:
\[
V_{r,\m}
=\{v\in V_r:\m^qv=0\text{ for some }q\ge1\}.
\]
Here $\m^qv=0$ means that every element of the ideal $\m^q$
annihilates $v$.
There are only finitely many occurring maximal ideals.
Each already occurs in $V_1$, and $V_{r,\m}$ is nonzero for all
$r\ge1$ (Proposition \ref{prop:density}).
After a finite extension of the coefficient field, which we again
denote by $\F$, we may assume that all occurring Hecke eigensystems
have eigenvalues in $\F$.
For each occurring $\m$, define
\begin{equation}\label{eq:Adef}
A_\m=\varprojlim_{r\ge1}
\im\bigl(\T\longrightarrow\End_\F(V_{r,\m})\bigr).
\end{equation}
This is a local algebra with maximal ideal $\n_\m$ and residue field
$A_\m/\n_\m\simeq\F$.
To study $A_\m$ using the $G$-action, we work with the smooth dual
\[
\Pi=(\Ht/\varpi)^\vee
=\Hom_\F^{\mathrm{cts}}(\Ht/\varpi,\F),
\]
where $\F$ has the discrete topology.
Write $\pi(g)$ for the action of $g\in G$ on $\Pi$.
The induced Hecke operators commute with every $\pi(g)$.

The level spaces defining $A_\m$ can be recovered as character
spaces in $\Pi$.
Let $P\subseteq\Iw$ be the subgroup of lower triangular matrices
and, for all $r\ge0$, put
\[
J_r=\pmat{\Z_p^\times}{p^r\Z_p}{p\Z_p}{\Z_p^\times}.
\]
For a subgroup $U\subseteq\Iw$, write $\Pi^{(U,\varepsilon)}$
for the vectors on which $U$ acts through $\varepsilon$.
Then there are $\T$-equivariant identifications
\[
\Sbar\simeq\Pi^{(P,\varepsilon)},\qquad
V_r\simeq\Pi^{(J_r,\varepsilon)}.
\]
These identifications allow us to use conjugation inside $G$
to study the action on the level spaces.

We distinguish the Hecke eigenspace associated with $\m$
\[
\Pi[\m]=\{v\in\Pi:\m v=0\}
\]
from its generalized Hecke eigenspace
\[
\Pi_\m=\bigcup_{q\ge1}\Pi[\m^q].
\]
At each level, they give
\[
V_{r,\m}=\Pi_\m^{(J_r,\varepsilon)},\qquad
V_{r,\m}[\m]=\Pi[\m]^{(J_r,\varepsilon)}.
\]

\subsection{Main results}\label{sec:obstacles}

Our first theorem bounds the power of the maximal ideal needed
to annihilate each finite level space.

\begin{theorem}\label{thm:A}
Let $p$ be odd, and let $\Ht$ be a projective augmented module with
a commuting $\T$-action on $\Ht/\varpi$.
For every occurring $\m$, there is an integer
$1\le b_\m\le\dim_\F(\Pi^{K_1})_\m$, where
$K_1=1+pM_2(\Z_p)$, such that, for all $r\ge0$,
\begin{equation}\label{eq:A1}
\n_\m^{\,\lceil2b_\m p^{r/2}\rceil}V_{r,\m}=0.
\end{equation}
\end{theorem}

The square-root bound was motivated by computations of Hecke
operators for the Hamilton quaternion algebra at $p=5$, using
the explicit formulas in \cite{Tr23}.

Since $\dim_\F V_{r,\m}=\rho_\m p^r+O(1)$ for some $\rho_\m>0$,
Theorem \ref{thm:A} gives a square-root bound for the nilpotence
exponent in terms of the dimension of the level space.
This is an analogue of Medvedovsky's nilpotence bound \cite{Me15}.

In our proof, we use the action of $\GL_2(\Q_p)$ on $\Pi$ to
conjugate the level subgroups and reduce the estimate to a
calculation in a finite group algebra.

We apply this nilpotence estimate to study the structure of the
completed local Hecke algebra.
For this purpose, we impose the additional hypothesis that the
Hecke eigenspace $\Pi[\m]$ has finite length as a
smooth $G$-representation.

\begin{theorem}\label{thm:B}
Let $p\ge5$, let $\Ht$ be a projective augmented module with a
commuting $\T$-action on $\Ht/\varpi$, and let $\m$ occur.
Suppose that $\Pi[\m]$ has finite length as a smooth
$G$-representation.
Then
\begin{enumerate}
\item
$A_\m$ is noetherian of Krull dimension at least two.
\item
There is a constant $c_\m>0$ such that
$\dim_\F(A_\m/\n_\m^q)\ge c_\m q^2$
for all sufficiently large integers $q$.
\item
If, in addition, $\dim_\F(\n_\m/\n_\m^2)\le2$, then
$A_\m\simeq\F[[X,Y]]$.
\end{enumerate}
\end{theorem}

In our abstract setting, we do not have $q$-expansions, so the
classical perfect pairing defined using the first Fourier
coefficient is unavailable.
Instead, we compare the dimension of each level space with
that of its Hecke eigenspace.

For $r\ge1$, if $\n_\m^qV_{r,\m}=0$, the socle inequality gives
\[
\dim_\F(A_\m/\n_\m^q)
\ge
\frac{\dim_\F V_{r,\m}}
{\dim_\F V_{r,\m}[\n_\m]}.
\]
We apply this inequality to the finite-dimensional spaces
$V_{r,\m}$.
The argument therefore requires no finite-dimensionality assumption
on the full Hecke eigenspace $\Sbar[\m]$.

Under the finite-length hypothesis, Hu's estimate
\cite[Theorem 4.1]{Hu21} gives
\[
\dim_\F V_{r,\m}[\n_\m]=O(r+1).
\]
Put $q_r=\lceil2b_\m p^{r/2}\rceil$.
Theorem \ref{thm:A} gives $\n_\m^{q_r}V_{r,\m}=0$, while
the dimension formula \eqref{eq:level-count} gives
\[
\dim_\F V_{r,\m}=\rho_\m p^r+O(1),
\qquad \rho_\m>0.
\]
The socle inequality therefore gives a constant $c>0$
such that, for all sufficiently large $r$,
\[
\frac{\dim_\F(A_\m/\n_\m^{q_r})}{q_r}
\ge c\frac{p^{r/2}}{r+1}\longrightarrow\infty.
\]
We prove noetherianness separately in Proposition
\ref{prop:FLnoeth}, using the projective-envelope formalism
of \cite{PT21}.
Hilbert--Samuel theory then gives the dimension bound and the
quadratic lower bound in Theorem \ref{thm:B}.
The additional cotangent bound yields the power series presentation.

We apply Theorem \ref{thm:B} in two arithmetic settings.
Fix a tame level $K^p$ and a non-Eisenstein maximal ideal $\m$.

For modular curves, $\Pi$ is the localized mod $p$ cohomology
of the tower with varying level at $p$:
\[
\Pi=\varinjlim_{K_p}
H^1_{\mathrm{\acute et}}
\bigl(Y(K^pK_p)_{\overline\Q},\F\bigr)_\m,
\]
where $Y(K^pK_p)$ is the modular curve of level $K^pK_p$ and
$K_p$ runs through compact open subgroups of $G$.

For a definite quaternion algebra $D/\Q$ split at $p$, we have
\[
\Pi=C^\infty\!\left(
D^\times\backslash(D\otimes_\Q\A_f)^\times/K^p,\F
\right),
\]
the space of locally constant $\F$-valued functions on this
quotient, where $\A_f$ denotes the finite adeles of $\Q$.
These are quaternionic automorphic forms of fixed tame level
and arbitrary level at $p$; the splitting of $D$ at $p$
gives the smooth $G$-action by right translation.

In both cases, Theorem \ref{thm:B} applies when the simultaneous
Hecke eigenspace $\Pi[\m]$ has finite length.
To state a sufficient condition in the modular case, let
$\bar\rho_\m:G_\Q\to\GL_2(\F)$ be the absolutely irreducible
residual Galois representation attached to the Hecke eigensystem
$\m$.
Emerton's mod $p$ local--global compatibility \cite{Em11}
gives finite length under the local condition on
$\bar\rho_\m|_{G_{\Q_p}}$ stated below.
For definite quaternion algebras, we assume finite length directly.

The power series conclusion requires two further conditions.
Let $\Gamma$ be the $p$-Sylow subgroup of the finite central
quotient determined by $K^p$ and $\varepsilon$
(Proposition \ref{prop:upper}); it measures the possible
$p$-power variation of the central character.
We call $\bar\rho_\m$ \emph{unobstructed} if
\[
H^2(G_{\Q,\mathcal R},\ad^0\bar\rho_\m)=0.
\]
Here $\mathcal R$ is the set of allowed ramification primes,
including $p$ and the bad primes of the arithmetic level,
$G_{\Q,\mathcal R}$ is the Galois group unramified outside
$\mathcal R$, and $\ad^0$ denotes the trace-zero adjoint
representation.
When $\Gamma=1$ and $\bar\rho_\m$ is unobstructed,
Corollaries \ref{cor:arithmetic} and \ref{cor:modular-arithmetic}
give the power series conclusion in the two settings.
We therefore obtain the following corollary.

\begin{corollary}
Let $p\ge5$, let $\m$ be a non-Eisenstein maximal ideal,
and choose $\varepsilon$ with $\Sbar^{(\varepsilon)}_\m\ne0$.
Suppose that we are in one of the following settings.
\begin{enumerate}
\item
\textup{(Modular curves.)}
Assume that, after extending scalars to $\overline\F$,
$\bar\rho_\m|_{G_{\Q_p}}$ is not isomorphic to any representation
of either of the forms
\[
\chi\otimes\begin{pmatrix}1&*\\0&1\end{pmatrix}
\qquad\text{or}\qquad
\chi\otimes\begin{pmatrix}1&*\\0&\omega\end{pmatrix},
\]
where $\chi:G_{\Q_p}\to\overline\F^\times$ is any continuous
character and $\omega$ is the mod $p$ cyclotomic character.
In both forms, $*$ may be zero.
\item
\textup{(Definite quaternion algebras.)}
Assume that the quaternion algebra is split at $p$, the tame
level is neat, and $\Pi[\m]$ has finite length as a smooth
$G$-representation.
\end{enumerate}
Then $A_\m$ is noetherian of Krull dimension at least two.
If, moreover, $\Gamma=1$ and $\bar\rho_\m$ is unobstructed, then
\[
A_\m\simeq\F[[X,Y]].
\]
\end{corollary}

\subsection{Notation and conventions}\label{sec:conventions}

Throughout, $p$ is an odd prime.
Let $E/\Q_p$ be finite, with ring of integers $\cO$,
uniformizer $\varpi$, and residue field $\F$.
We enlarge the coefficient fields when necessary,
extend modules and Hecke actions by scalars,
and retain the same notation.

Compact group-ring modules are right modules, and their
duals are left representations.
Maps from compact modules to smooth representations
are understood to be continuous.
Hecke algebras are commutative and unital.

For a local ring $(B,\n)$ and a $B$-module $L$, write
\[
L[\n]=\{v\in L:\n v=0\}
\]
for its socle.
We use $w$ for the weight variable, so that
$\cO[[1+p\Z_p]]\simeq\cO[[w]]$.

For a finite set of primes $\mathcal R$ containing $p$,
let $G_{\Q,\mathcal R}$ denote the Galois group of the
maximal extension of $\Q$ unramified at finite primes
outside $\mathcal R$.
For $\ell\notin\mathcal R$, we use arithmetic Frobenius,
so $\chi_{\mathrm{cyc}}(\Frob_\ell)=\ell$, and normalize
Galois representations by
\[
\det\bigl(X-\rho(\Frob_\ell)\bigr)
=X^2-T_\ell X+\ell S_\ell.
\]
We write $\omega$ for the mod $p$ cyclotomic character,
$\bar\rho_p=\bar\rho|_{G_{\Q_p}}$, and
$\mathrm{Teich}$ for the Teichm\"uller lift.

\section{Abstract automorphic forms and their filtrations}
\label{sec:setup}

We define the space $\Sbar$ of residual automorphic forms and its
level spaces $V_r$. We identify $\Sbar$ and $V_r$ with the character
spaces $\Pi^{(P,\varepsilon)}$ and $\Pi^{(J_r,\varepsilon)}$,
respectively. The first description gives the Mahler filtration, whose successive
quotients allow us to compute $\dim_\F V_{r,\m}$. The second lets us
conjugate $J_r$ inside $G$ in the proof of Theorem \ref{thm:A}.
We finish with examples from modular curves and definite quaternion
algebras.

\subsection{The space of residual automorphic forms}
\label{sec:residual}
Put $G=\GL_2(\Q_p)$, $K=\GL_2(\Z_p)$, and
$K_n=1+p^nM_2(\Z_p)$ for $n\ge1$.
Define
\[
\Iw=\pmat{\Z_p^\times}{\Z_p}{p\Z_p}{\Z_p^\times},\qquad
\Iwone=\pmat{1+p\Z_p}{\Z_p}{p\Z_p}{1+p\Z_p},\qquad
P=\pmat{\Z_p^\times}{0}{p\Z_p}{\Z_p^\times}.
\]
These are the Iwahori subgroup, its pro-$p$ Sylow subgroup, and its
lower triangular part. The group $P$ is denoted
$B^{\mathrm{op}}(\Z_p)$ in \cite[(2.4.1)]{LTXZ26}.
Let $T,N,N^-$ be the diagonal torus and the upper and lower
unipotent subgroups, and write
\[
n(z)=\smat1z01,\qquad n^-(c)=\smat10c1,\qquad
\alpha=\smat p001.
\]

Consider the surjective homomorphism
\[
\Iw\longrightarrow T(\F_p)\simeq
\F_p^\times\times\F_p^\times,
\qquad
\smat abcd\longmapsto(\bar a,\bar d).
\]

Fix a character $\varepsilon:T(\F_p)\to\F^\times$.
We regard it as a character of $\Iw$ by setting
\[
\varepsilon\!\left(\smat abcd\right)
=\varepsilon(\bar a,\bar d),
\]
and use the same notation for its restriction to any subgroup
of $\Iw$.

For each $n\in\Z$, define a character $\chi_n:\Iw\to\F^\times$ by
\[
\chi_n\!\left(\smat abcd\right)
=\left(\frac{\bar a}{\bar d}\right)^n.
\]
Since $\bar a/\bar d\in\F_p^\times$, the character $\chi_n$
depends only on $n$ modulo $p-1$.

\begin{definition}
A \emph{projective augmented module} is a compact right
$\cO[[\Iw]]$-module $\Ht$, finitely generated and projective over
$\cO[[\Iw]]$, with a continuous right $G$-action by $\cO$-linear
maps extending the $\Iw$-action.
\end{definition}

The $\cO[[K]]$-projective augmented modules of
\cite[Definition 2.2]{LTXZ26} satisfy our definition.
Indeed, $\cO[[K]]$ is free of rank $p+1$ over
$\cO[[\Iw]]$, so restriction from $K$ to $\Iw$
preserves finite generation and projectivity.

We write $\Ht/\varpi$ for $\Ht/\varpi\Ht$.
This is a finitely generated projective $\F[[\Iw]]$-module.
Its restriction to $\F[[\Iwone]]$ is free of finite rank,
since $\Iwone$ is pro-$p$ and $\F[[\Iwone]]$ is local.

Set
\[
\Pi=\Hom^{\mathrm{cont}}_\F(\Ht/\varpi,\F).
\]
The right $G$-action on $\Ht/\varpi$ induces a left
$G$-action on $\Pi$, given by
\[
(\pi(g)\lambda)(x)=\lambda(xg),
\qquad
g\in G,\quad \lambda\in\Pi,\quad x\in\Ht/\varpi.
\]

For a closed subgroup $U\subseteq\Iw$ and a continuous
character $\eta:U\to\F^\times$, define
\[
\Pi^{(U,\eta)}
=\{F\in\Pi:\pi(u)F=\eta(u)F
\text{ for all }u\in U\}.
\]
We also consider the induced representation
\[
\Ind_U^{\Iw}\eta
=
\left\{
f\in C^0(\Iw,\F):
f(gu)=\eta(u)f(g)
\text{ for all }g\in\Iw,\ u\in U
\right\},
\]
with right $\Iw$-action given by
$(f\cdot g')(g)=f(g'g)$.

The following lemma records the properties of $\Pi$
that we will use.

\begin{lemma}\label{lem:Pi}
\begin{enumerate}
\item
$\Pi$ is smooth and admissible.
\item
For every closed subgroup $U\subseteq\Iw$ and continuous
character $\eta:U\to\F^\times$, there are natural
isomorphisms
\[
\Pi^{(U,\eta)}
\simeq\Hom_{\F[[U]]}(\Ht/\varpi,\eta)
\simeq\Hom_{\F[[\Iw]]}
(\Ht/\varpi,\Ind_U^{\Iw}\eta).
\]
\item
The functor $\Hom_{\F[[\Iw]]}(\Ht/\varpi,-)$
is exact on smooth right $\Iw$-representations.
\end{enumerate}
\end{lemma}

\begin{proof}
For (1), by finite generation and continuity,
every $F\in\Pi$ factors through $(\Ht/\varpi)_V$
for some open normal subgroup $V\subseteq\Iw$.
Hence $\Pi= \bigcup_{V\triangleleft_{\mathrm{open}}\Iw}\Pi^V,$
so $\Pi$ is smooth.

For any open $V\subseteq\Iw$, we have
$\Pi^V\simeq\Hom_\F((\Ht/\varpi)_V,\F)$ and
$\dim_\F(\Ht/\varpi)_V<\infty$,
since $\Ht/\varpi$ is finitely generated over $\F[[V]]$.

For any compact open subgroup $W\subseteq G$,
the intersection $W\cap\Iw$ is open in $\Iw$, and
$\Pi^W\subseteq\Pi^{W\cap\Iw}.$
Thus $\dim_\F\Pi^W<\infty$, proving that $\Pi$
is admissible.

For (2), the first identification follows from the
dual action. The second is given by
\[
F\longmapsto\phi_F,\qquad
\phi_F(x)(g)=F(xg),
\]
with inverse given by evaluation at the identity.
Indeed,
\[
\phi_F(x)(gu)=\eta(u)\phi_F(x)(g),
\qquad
\phi_F(xh)=\phi_F(x)\cdot h.
\]
Choose an open normal subgroup $V\subseteq\Iw$
fixing $F$. Then
\[
\phi_F(x)(vg)=\phi_F(x)(g)
\qquad (v\in V).
\]
Thus $\phi_F(x)$ is determined by finitely many
continuous functions $x\mapsto F(xg)$, proving
continuity of $\phi_F$.
This applies to any closed $U\subseteq\Iw$.

For (3), projectivity gives
\[
\F[[\Iw]]^d\simeq(\Ht/\varpi)\oplus Q
\]
for some $\F[[\Iw]]$-module $Q$.
Since
\[
\Hom_{\F[[\Iw]]}(\F[[\Iw]]^d,M)\simeq M^d
\]
for every smooth right $\Iw$-representation $M$,
the asserted exactness follows.
\end{proof}

We now identify the induced representation
$\Ind_P^{\Iw}\varepsilon$ with $C^0(\Z_p,\F)$
and compute the corresponding Iwahori action.

Every element of $\Iw$ can be written uniquely as
$n(z)v$, with $z\in\Z_p$ and $v\in P$.
Thus the decomposition $\Iw=N(\Z_p)P$ identifies
$\Iw/P$ with $\Z_p$.
Restricting functions to $n(\Z_p)$ gives
\[
\Ind_P^{\Iw}\varepsilon\simeq C^0(\Z_p,\F).
\]

We now compute the Iwahori action on $C^0(\Z_p,\F)$.
Let $f\in\Ind_P^{\Iw}\varepsilon$. Then
\[
f(gv)=\varepsilon(v)f(g),
\qquad v\in P.
\]
The corresponding function $h\in C^0(\Z_p,\F)$
is given by $h(z)=f(n(z))$.

For $u=\smat abcd\in\Iw$, we have
\[
u\,n(z)
=
n\Bigl(\frac{az+b}{cz+d}\Bigr)
\smat{(ad-bc)/(cz+d)}{0}{c}{cz+d}.
\]
Since $cz+d\in\Z_p^\times$, the second factor
lies in $P$, and its diagonal entries reduce to
$\bar a$ and $\bar d$.
Consequently,
\begin{equation}\label{eq:action}
\begin{aligned}
(h\|_u)(z)
&=f(un(z))
&=\varepsilon(\bar a,\bar d)
h\Bigl(\frac{az+b}{cz+d}\Bigr).
\end{aligned}
\end{equation}

We can now define the residual space.

\begin{definition}
The \emph{residual space} is
\[
\Sbar=\Sbar^{(\varepsilon)}(\Ht)
:=\Hom_{\F[[\Iw]]}
(\Ht/\varpi,\Ind_P^{\Iw}\varepsilon).
\]
\end{definition}

By Lemma \ref{lem:Pi}(2), there is a natural identification
$\Sbar\simeq\Pi^{(P,\varepsilon)}.$
Explicitly, an element $F\in\Pi^{(P,\varepsilon)}$
corresponds to the residual form sending
$x\in\Ht/\varpi$ to the function
\[
z\longmapsto F(xn(z)).
\]
Thus the descriptions using coefficient functions
and smooth representations agree.

We next relate this space to the abstract
$p$-adic automorphic forms of \cite{LTXZ26}.
By projectivity of $\Ht$, specialization commutes
with the $\Hom$ construction defining
$S^{(\varepsilon)}_{p\text{-adic}}(\Ht)$.
Consequently,
\[
S^{(\varepsilon)}_{p\text{-adic}}(\Ht)/(\varpi,w)
\simeq\Sbar,
\]
where the universal character specializes to
$\varepsilon$
\cite[(2.4.2)--(2.4.4)]{LTXZ26}.

Finally, let $\T$ be a commutative $\F$-algebra of continuous
$G$-equivariant endomorphisms of $\Ht/\varpi$.
It acts on $\Pi$ and $\Sbar$ by precomposition:
\[
(t\phi)(x)=\phi(tx).
\]
Since the Hecke action commutes with $G$, it preserves
the character spaces.
In the arithmetic examples, $\T$ is generated by the
operators $T_\ell$ and $S_\ell$ for
$\ell\notin\mathcal R$, where $p\in\mathcal R$.
The abstract arguments use only commutativity and
$G$-equivariance.

\subsection{The level and Mahler filtrations}
\label{sec:level}
\label{sec:mahler}\label{sec:modified}\label{sec:graded}

The level filtration describes forms of increasing $p$-power
level. We refine it using the degree of coefficient functions
and identify the successive quotients with Iwahori character
spaces.

For $r\ge0$ and $m\ge1$, put
\[
J_r=\pmat{\Z_p^\times}{p^r\Z_p}{p\Z_p}{\Z_p^\times},
\qquad
\Gamma_0(p^m)
=\pmat{\Z_p^\times}{\Z_p}{p^m\Z_p}{\Z_p^\times}.
\]

Then $J_0=\Iw$, $J_r\supseteq K_r$ for $r\ge1$,
$\bigcap_rJ_r=P$, and
\begin{equation}\label{eq:Jr-conj}
J_r=\alpha^r\Gamma_0(p^{r+1})\alpha^{-r}.
\end{equation}
Since $J_r=N(p^r\Z_p)P$, the coset space $\Iw/J_r$
identifies with $\Z/p^r\Z$ by
$n(z)J_r\mapsto z\bmod p^r$.

\begin{definition}
For every $r\ge0$, put
\[
V_r=
\Hom_{\F[[\Iw]]}
\bigl(\Ht/\varpi,\Map(\Z/p^r\Z,\F)\bigr)
\subseteq\Sbar,
\]
where the coefficient action is given by \eqref{eq:action}.
The spaces $V_r$ form the \emph{level filtration}.
\end{definition}

\begin{lemma}\label{lem:level}
\begin{enumerate}
\item
For every $r\ge0$, there is an isomorphism of
right $\Iw$-representations
\[
\Ind_{J_r}^{\Iw}\varepsilon
\xrightarrow{\sim}\Map(\Z/p^r\Z,\F),
\qquad
f\longmapsto\bigl(z\longmapsto f(n(z))\bigr).
\]
\item
For every $r\ge0$, evaluation of coefficient functions at $0$
induces a natural $\T$-equivariant isomorphism
\[
V_r\xrightarrow{\sim}\Pi^{(J_r,\varepsilon)},
\qquad
\phi\longmapsto\bigl(x\longmapsto\phi(x)(0)\bigr).
\]
\item
Each $V_r$ is finite-dimensional and $\T$-stable, and
$\Sbar=\bigcup_{r\ge0}V_r.$
\end{enumerate}
\end{lemma}

\begin{proof}
For (1), restriction to $n(\Z_p)$ gives the isomorphism via
$\Iw/J_r\simeq\Z/p^r\Z$, since $\varepsilon$ is trivial on
$N(p^r\Z_p)$. The action is \eqref{eq:action}.

Part (2) follows from (1) and Lemma \ref{lem:Pi}(2).
Evaluation at the identity corresponds to evaluation at $0$.

For (3), since $K_{r+1}\subseteq J_r$ and
$\varepsilon|_{K_{r+1}}=1$, admissibility gives
\[
\dim_\F V_r
=\dim_\F\Pi^{(J_r,\varepsilon)}
\le\dim_\F\Pi^{K_{r+1}}<\infty.
\]
The $\T$-stability follows because $\T$ commutes with $G$.

Finally, smoothness and $PK_n=J_n$ give
\[
\Pi^{(P,\varepsilon)}
=\bigcup_{n\ge1}
\bigl(\Pi^{(P,\varepsilon)}\cap\Pi^{K_n}\bigr)
=\bigcup_{n\ge1}\Pi^{(J_n,\varepsilon)}.
\]
Using $\Sbar\simeq\Pi^{(P,\varepsilon)}$ and part~(2), we obtain
$\Sbar=\bigcup_{r\ge0}V_r.$
\end{proof}

We use the usual Mahler basis to describe the coefficient space
at level $r$. The reductions of its first $p^r$ elements form
a basis of $\Map(\Z/p^r\Z,\F)$. The modified Mahler basis gives
the same degree filtration and lets us compute the Iwahori action
on its successive quotients.

Let $b_n(z)=\binom zn$ be the usual Mahler basis.
Following \cite[(3.14.1)--(3.14.2)]{LTXZ26},
define the modified Mahler polynomials $m_n(z)$ by
\[
f_0(z)=z,\qquad
f_{i+1}(z)=\frac{f_i(z)^p-f_i(z)}p,\qquad
m_n(z)=\prod_{i\ge0}f_i(z)^{n_i},
\]
where $n=\sum_i n_ip^i$ and $0\le n_i<p$.
Each $f_i$ is an integer-valued polynomial of degree $p^i$.

\begin{proposition}\label{prop:mahler-bases}
\begin{enumerate}
\item
For every $N\ge0$,
\[
\bigoplus_{n=0}^N\Z_pm_n
=
\bigoplus_{n=0}^N\Z_pb_n
=
\{h\in\Q_p[z]:
\deg h\le N,\ h(\Z_p)\subseteq\Z_p\}.
\]
\item
For every $r\ge0$, the reductions modulo $p$
of $b_0,\ldots,b_{p^r-1}$ and
$m_0,\ldots,m_{p^r-1}$, extended to $\F$,
form bases of $\Map(\Z/p^r\Z,\F)$.
\end{enumerate}
\end{proposition}

\begin{proof}
We first prove (1).
Let $\Delta h(z)=h(z+1)-h(z)$ be the forward
difference operator.
Every polynomial $h$ of degree at most $N$ has the expansion
\[
h(z)=\sum_{n=0}^N(\Delta^nh)(0)b_n(z).
\]
If $h$ is integer-valued on $\Z_p$, then
$(\Delta^nh)(0)\in\Z_p$.
Thus $b_0,\ldots,b_N$ form a $\Z_p$-basis of the
integer-valued polynomials of degree at most $N$.

The polynomials $f_i$ are integer-valued by induction,
since $x^p-x\in p\Z_p$ for every $x\in\Z_p$.
Their leading coefficients are
$p^{-(1+p+\cdots+p^{i-1})}$.
Hence, writing $n=\sum_i n_ip^i$, the leading
coefficient of $m_n$ is
\[
p^{-\sum_i n_i(1+p+\cdots+p^{i-1})}
=p^{-v_p(n!)}
\]
by Legendre's formula.
Since $b_n$ has leading coefficient $1/n!$,
the change of basis from $(m_n)$ to $(b_n)$
is triangular with diagonal entries
$n!p^{-v_p(n!)}\in\Z_p^\times$.
This proves (1).

For (2), let $0\le a,j<p$ and $k\ge0$.
Comparing coefficients in
$(1+X)^{j+pz}\equiv(1+X)^j(1+X^p)^z\bmod p$
gives
\[
b_{pk+a}(j+pz)
\equiv\binom ja\,b_k(z)\pmod p.
\]
One may first make this coefficient comparison for nonnegative
integers $z$ and then extend it to $\Z_p$ by continuity.
Iterating this identity shows that $b_n\bmod p$
depends only on $z\bmod p^r$ whenever $n<p^r$.
Moreover, the evaluation matrix
$(b_n(j))_{0\le j,n<p^r}$ is triangular
with diagonal entries $1$.
Thus these $p^r$ functions form a basis of
$\Map(\Z/p^r\Z,\F)$.
The same holds for $(m_n)$ by the triangular
change of basis.
\end{proof}

We use $m_n$ and $b_n$ also for their reductions modulo $p$,
extended to $\F$. Define the degree filtration on $C^0(\Z_p,\F)$ by
\[
C^0_{\deg\le N}
=
\Span_\F\{m_0,\ldots,m_N\}
=
\Span_\F\{b_0,\ldots,b_N\},
\]
with $C^0_{\deg\le-1}=0$.
For a nonzero function $h$, its \emph{Mahler degree}
is the least $N$ such that $h\in C^0_{\deg\le N}$.
In particular,
\[
C^0_{\deg\le p^r-1}
=\Map(\Z/p^r\Z,\F).
\]

We now show that the Iwahori action preserves the
Mahler filtration and determine its action on the
successive quotients.

Let $\bar f_i:\Z_p\to\F_p$ denote the reduction of
$f_i$ modulo $p$.
The following lemma describes how these functions transform
under the Iwahori action.

\begin{lemma}\label{lem:digits}
\begin{enumerate}
\item
For every $r\ge1$, the map
$z\bmod p^r
\longmapsto
(\bar f_0(z),\ldots,\bar f_{r-1}(z))$
is a bijection from $\Z/p^r\Z$ onto $\F_p^r$.
\item
For $g=\smat abcd\in\Iw$ and every $i\ge0$, we have
\begin{equation}\label{eq:digit-action}
\bar f_i\Bigl(\frac{az+b}{cz+d}\Bigr)
=
\frac{\bar a}{\bar d}\bar f_i(z)
+
H_i(\bar f_0(z),\ldots,\bar f_{i-1}(z)),
\end{equation}
where $H_i:\F_p^i\to\F_p$ depends on $g$,
and $H_0=\bar b/\bar d$.
\end{enumerate}
\end{lemma}

\begin{proof}
For (1), the binomial theorem gives, for $z,t\in\Z_p$,
\[
f_1(z+p^jt)-f_1(z)
\equiv-p^{j-1}t\pmod{p^j}
\]
for $j\ge1$.
Iterating, for $j\ge i$ we obtain
\[
f_i(z+p^jt)-f_i(z)
\equiv(-1)^ip^{j-i}t\pmod{p^{j-i+1}}.
\]
Thus $\bar f_i$ depends only on $z\bmod p^{i+1}$,
and
\[
\bar f_i(z+p^it)
=
\bar f_i(z)+(-1)^i\bar t.
\]
Hence $\bar f_i$ distinguishes the $p$ possible lifts
of each residue class modulo $p^i$, proving (1)
by induction.

For (2), write $g(z)=(az+b)/(cz+d)$.
For $i\ge1$, we have
\[
g(z+p^it)-g(z)
=
\frac{(ad-bc)p^it}
{(cz+d)(c(z+p^it)+d)}
\equiv
\frac{a}{d}p^it\pmod{p^{i+1}}.
\]
By (1),
\[
\bar f_i(g(z+p^it))-\bar f_i(g(z))
=
(-1)^i\frac{\bar a}{\bar d}\bar t.
\]
Therefore
$\bar f_i(g(z))-(\bar a/\bar d)\bar f_i(z)$
depends only on $z\bmod p^i$.
By (1), this is equivalent to depending only on
$\bar f_0(z),\ldots,\bar f_{i-1}(z)$.
For $i=0$, the assertion follows from
\[
g(z)
\equiv
(\bar a/\bar d)z+\bar b/\bar d\pmod p.
\]
\end{proof}

The following proposition shows that the Iwahori action
preserves Mahler degree.

\begin{proposition}\label{prop:mahler-action}
For every $n\ge0$, the following hold:
\begin{enumerate}
\item
For every $u\in\Iw$,
\[
m_n\|_u-\varepsilon(u)\chi_n(u)m_n
\in C^0_{\deg\le n-1}.
\]
\item
The spaces $C^0_{\deg\le N}$ are $\Iw$-stable
for all $N\ge0$.
\item
There is an isomorphism of $\Iw$-representations
\[
C^0_{\deg\le n}/C^0_{\deg\le n-1}
\simeq\varepsilon\chi_n.
\]
\end{enumerate}
\end{proposition}

\begin{proof}
For (1), write $n=\sum_i n_ip^i$ and
$u=\smat abcd\in\Iw$.
By \eqref{eq:digit-action},
\[
m_n\|_u
=
\varepsilon(u)\prod_i
\left(
\frac{\bar a}{\bar d}\bar f_i
+
H_i(\bar f_0,\ldots,\bar f_{i-1})
\right)^{n_i}.
\]
The leading term is
\[
\varepsilon(u)
\left(\frac{\bar a}{\bar d}\right)^{\sum_i n_i}m_n
=
\varepsilon(u)\chi_n(u)m_n,
\]
since $n\equiv\sum_i n_i\pmod{p-1}$.

We claim that every other term has Mahler degree less
than $n$.
Let $i$ be the largest index at which an $H_i$ occurs.
Such a term is, up to a scalar, of the form
\[
\left(\prod_{j>i}\bar f_j^{n_j}\right)
\bar f_i^k
Q(\bar f_0,\ldots,\bar f_{i-1}),
\qquad 0\le k<n_i,
\]
where $Q$ is a function of the lower coordinates.
By Lemma \ref{lem:digits}(1) and polynomial interpolation
on $\F_p^i$, this function is a linear combination
of monomials
\[
\prod_{j<i}\bar f_j^{e_j},
\qquad 0\le e_j<p.
\]
Each resulting monomial therefore has Mahler degree
at most
\[
\sum_{j>i}n_jp^j+kp^i+(p^i-1)
\le
\sum_{j>i}n_jp^j+n_ip^i-1
\le n-1.
\]
This proves (1).

Statement (2) follows immediately from (1).
Finally, the quotient of $C^0_{\deg\le n}$ by
$C^0_{\deg\le n-1}$ is one-dimensional, generated by the class of $m_n$.
By (1), $\Iw$ acts on this class through
$\varepsilon\chi_n$, proving (3).
\end{proof}

We now define a filtration on residual automorphic forms.

\begin{definition}
The \emph{Mahler filtration} on $\Sbar$ is
\[
\Sbar_{\deg\le N}
=
\Hom_{\F[[\Iw]]}
(\Ht/\varpi,C^0_{\deg\le N}).
\]
We set $\Sbar_{\deg\le-1}=0$.
\end{definition}

The Mahler filtration is $\T$-stable, and
$V_r=\Sbar_{\deg\le p^r-1}.$
Since $\Sbar=\bigcup_{r\ge0}V_r$, every residual form lies in
$\Sbar_{\deg\le N}$ for some $N$.

\begin{proposition}\label{prop:graded}
For every $n\ge0$, let $0\le a\le p-2$ satisfy
$a\equiv n\pmod{p-1}$.
There is a natural $\T$-equivariant isomorphism
\[
\Sbar_{\deg\le n}/\Sbar_{\deg\le n-1}
\simeq\Pi^{(\Iw,\varepsilon\chi_a)}.
\]
Consequently,
\[
\dim_\F\Sbar_{\deg\le N}
=
\sum_{n=0}^N
\dim_\F\Pi^{(\Iw,\varepsilon\chi_n)}.
\]
\end{proposition}

\begin{proof}
By Proposition \ref{prop:mahler-action}(3),
we have an exact sequence of $\Iw$-representations
\[
0\longrightarrow C^0_{\deg\le n-1}
\longrightarrow C^0_{\deg\le n}
\longrightarrow\varepsilon\chi_n
\longrightarrow0.
\]
Applying the exact functor
$\Hom_{\F[[\Iw]]}(\Ht/\varpi,-)$
from Lemma \ref{lem:Pi}(3) gives
\[
\Sbar_{\deg\le n}/\Sbar_{\deg\le n-1}
\simeq
\Hom_{\F[[\Iw]]}(\Ht/\varpi,\varepsilon\chi_n).
\]
By Lemma \ref{lem:Pi}(2), the latter is
$\Pi^{(\Iw,\varepsilon\chi_a)}$,
since $\chi_n=\chi_a$.
These identifications are $\T$-equivariant.
Summing dimensions gives the formula.
\end{proof}

\subsection{Residual Hecke eigensystems and dimension growth}
\label{sec:primary}

We study the action of $\T$ on the space
$\Sbar=\Pi^{(P,\varepsilon)}$
of residual $p$-adic automorphic forms.
By Proposition \ref{prop:graded}, for every $n\ge0$,
\[
\Sbar_{\deg\le n}/\Sbar_{\deg\le n-1}
\simeq
\Pi^{(\Iw,\varepsilon\chi_n)}.
\]
We localize the degree filtration at the occurring maximal
ideals $\m$ of $\T$ and compute the dimensions of the local
components $\Sbar_{\deg\le N,\m}$ and $V_{r,\m}$.

Since $\chi_n$ depends only on $n\bmod(p-1)$, every residual
Hecke eigensystem occurring in $\Sbar$ occurs in the
finite-dimensional space
\[
\bigoplus_{a=0}^{p-2}\Pi^{(\Iw,\varepsilon\chi_a)}.
\]
Its Hecke image is a finite-dimensional $\F$-algebra, so we
may enlarge $\F$ finitely so that all these eigensystems
take values in $\F$.
We extend all spaces and Hecke actions by scalars and
retain the same notation.

Let $\Sigma$ be the finite set of maximal ideals of $\T$
that occur in $\Sbar$.
Our choice of $\F$ ensures that $\T/\m=\F$
for every $\m\in\Sigma$.

\begin{definition}
A \emph{residual Hecke eigensystem} is an $\F$-algebra
homomorphism
\[
\lambda_\m:\T\longrightarrow\F,
\qquad
\ker(\lambda_\m)=\m.
\]
We identify the eigensystem with its kernel.
For $\m\in\Sigma$, put
$\Sbar_\m
=
\bigcup_{q\ge1}\Sbar[\m^q].$
We call $\Sbar_\m$ the \emph{local component} at $\m$.
Proposition \ref{prop:density} identifies it canonically
with the localization of $\Sbar$ at $\m$.

The \emph{Hecke eigenspace} belonging to $\lambda_\m$ is
\[
\Sbar[\m]
=
\{f\in\Sbar:
Tf=\lambda_\m(T)f\text{ for every }T\in\T\}.
\]
Thus $\Sbar[\m]$ consists of vectors annihilated by $\m$,
while $\Sbar_\m$ also includes generalized Hecke
eigenvectors.
\end{definition}

We intersect $\Sbar_\m$ with the Mahler and level
filtrations:
\[
\Sbar_{\deg\le N,\m}
=
\Sbar_{\deg\le N}\cap\Sbar_\m,
\qquad
V_{r,\m}
=
V_r\cap\Sbar_\m.
\]
For $0\le a\le p-2$, define $d_{\m,a}
=
\dim_\F
\left(\Pi^{(\Iw,\varepsilon\chi_a)}\right)_\m.$
By Proposition \ref{prop:graded}, $d_{\m,a}$ is the
dimension contributed by $\m$ in each degree
$n\equiv a\pmod{p-1}$.

\begin{proposition}\label{prop:density}
\begin{enumerate}
\item
There is a decomposition
$\Sbar
=
\bigoplus_{\m\in\Sigma}\Sbar_\m,$
and the natural map
$
\Sbar_\m
\longrightarrow
\T_\m\otimes_\T\Sbar,
\quad
f\longmapsto1\otimes f,
$
is an isomorphism.
\item
For every $\m\in\Sigma$ and $N\ge0$, we have
\[
\dim_\F\Sbar_{\deg\le N,\m} =
(N+1)\rho_\m+O(1), \textrm{ where }\rho_\m =
\frac{1}{p-1}\sum_{a=0}^{p-2}d_{\m,a}>0.
\]
\item
The residual Hecke eigensystems occurring in $\Sbar$
are exactly those indexed by $\Sigma$.
Every one of them occurs in $V_r$ for every $r\ge1$.
\end{enumerate}
\end{proposition}

\begin{proof}
By Proposition \ref{prop:graded}, the space
$\Sbar_{\deg\le N}$ has a finite filtration with quotients
$\Pi^{(\Iw,\varepsilon\chi_n)}$ for $0\le n\le N$.
All its Hecke eigensystems are indexed by $\Sigma$, since
$\Sbar_{\deg\le N}\subseteq\Sbar$.

For (1), the Hecke image on $\Sbar_{\deg\le N}$ is a
finite-dimensional commutative $\F$-algebra.
Its decomposition into local factors gives
\[
\Sbar_{\deg\le N}
=\bigoplus_{\m\in\Sigma}\Sbar_{\deg\le N,\m},
\]
with each summand canonically identified with the localization
at $\m$.
The inclusions between successive filtration steps preserve
these summands.
Taking the union over $N$ gives the asserted decomposition.
Since localization commutes with this union, the natural map
$\Sbar_\m\to\T_\m\otimes_\T\Sbar$ is an isomorphism.

For (2), exactness of localization and
Proposition \ref{prop:graded} give
\[
0\longrightarrow\Sbar_{\deg\le n-1,\m}
\longrightarrow\Sbar_{\deg\le n,\m}
\longrightarrow
\left(\Pi^{(\Iw,\varepsilon\chi_n)}\right)_\m
\longrightarrow0.
\]
The last term has dimension $d_{\m,a}$ when
$n\equiv a\pmod{p-1}$.
Summing over $0\le n\le N$ gives
\[
\dim_\F\Sbar_{\deg\le N,\m}
=
\sum_{n=0}^{N}d_{\m,n\bmod(p-1)}.
\]
Each residue class modulo $p-1$ occurs
$(N+1)/(p-1)+O(1)$ times, proving the asserted asymptotic.

For $\m\in\Sigma$, part (1) gives $\Sbar_\m\ne0$.
Thus $\Sbar_{\deg\le N,\m}\ne0$ for some $N$.
The dimension sum above implies that at least one
$d_{\m,a}$ is positive, so $\rho_\m>0$.

For (3), the first assertion follows from the definition
of $\Sigma$.
Since $V_1=\Sbar_{\deg\le p-1}$, the dimension sum gives
\[
\dim_\F V_{1,\m}
=
2d_{\m,0}+\sum_{a=1}^{p-2}d_{\m,a}>0.
\]
A nonzero finite-dimensional generalized $\m$-eigenspace
contains a nonzero vector annihilated by $\m$.
Thus each eigensystem indexed by $\Sigma$ occurs in $V_1$,
and hence in every $V_r$ with $r\ge1$.
\end{proof}

Using $V_r=\Sbar_{\deg\le p^r-1}$ and the dimension sum
in the proof of Proposition \ref{prop:density}, we obtain
an exact formula for $\dim_\F V_{r,\m}$ at every level $r\ge0$.
Among $0,\ldots,p^r-1$, each residue class modulo $p-1$
occurs $(p^r-1)/(p-1)$ times, with one additional
occurrence of the class $0$.
Thus
\begin{equation}\label{eq:level-count}
\dim_\F V_{r,\m}
=
d_{\m,0}
+
\frac{p^r-1}{p-1}
\sum_{a=0}^{p-2}d_{\m,a}
=
\rho_\m p^r+d_{\m,0}-\rho_\m.
\end{equation}

We also compare the Hecke eigenspaces in the degree and
level filtrations.
For $r\ge0$ and $0\le N<p^r$, we have
\begin{equation}\label{eq:socle-growth}
\Sbar_{\deg\le N,\m}[\m]
=
\Sbar_{\deg\le N}\cap\Pi[\m]^{(P,\varepsilon)}
\subseteq
\Pi[\m]^{(J_r,\varepsilon)}.
\end{equation}
The equality follows from the identification
$\Sbar\simeq\Pi^{(P,\varepsilon)}$ and the fact that every
vector annihilated by $\m$ belongs to $\Sbar_\m$.
The inclusion follows from $\Sbar_{\deg\le N}\subseteq V_r$
and $V_r\simeq\Pi^{(J_r,\varepsilon)}$.
We use \eqref{eq:level-count} and \eqref{eq:socle-growth}
in the dimension argument in Section \ref{sec:duality}.

\subsection{Completed homology and its arithmetic realizations}
\label{sec:examples}
\label{sec:quaternion-ex}\label{sec:modular-ex}

\begin{example}[Definite quaternion algebras]\label{ex:quaternion}
Let $D/\Q$ be definite and split at $p$, with a fixed isomorphism
$D\otimes\Q_p\simeq M_2(\Q_p)$, and let
$K^p\subseteq(D\otimes\A_f^{(p)})^\times$ be compact open.
Write $D_f^\times=(D\otimes_\Q\A_f)^\times$ and
$D_p^\times=(D\otimes_\Q\Q_p)^\times$.
Put
\[
X=D^\times\backslash D_f^\times/K^p.
\]
Assume $K^p\Iw$ is neat in the sense of \cite[\S2.4]{LWX17}:
for every $x\in D_f^\times$,
\[
x^{-1}D^\times x\cap K^p\Iw=\{1\}.
\]
Then $\Iw$ acts freely on $X$, and $X/\Iw$ is finite,
say of cardinality $t$.
The module
\[
\Ht=\cO[[X]]=\varprojlim_n\cO[X/K_n]
\]
is free of rank $t$ over $\cO[[\Iw]]$, and
$\Pi=C^\infty(X,\F).$
Orbit representatives $\gamma_1,\ldots,\gamma_t$ give
\[
\Sbar\simeq C^0(\Z_p,\F)^{\oplus t},\qquad
\phi\longmapsto(\phi(\gamma_i))_i.
\]
Here $\gamma_i$ denotes also the Dirac class of that point in
$\Ht/\varpi$, so $\phi(\gamma_i)$ is a coefficient function on
$\Z_p$.
This is the residual version of \cite[(2.11.1)]{LWX17}.

If $K^pK$ is neat, $\Ht$ is also free over $\cO[[K]]$.
This is the degree-zero completed homology example of
\cite[Example 8.4]{LTXZ26}.

Let $\mathcal R$ contain $p$, the ramified primes of $D$,
and the primes where $K^p$ is not maximal.
For $\ell\notin\mathcal R$, choose
\[
\GL_2(\Z_\ell)\smat\ell001\GL_2(\Z_\ell)
=\bigsqcup_{\beta=1}^{\ell+1}\GL_2(\Z_\ell)v_\beta.
\]
Regard each $v_\beta$ as an element of $D_f^\times$ with
identity components away from $\ell$.
Set
\[
(T_\ell F)(x)=\sum_{\beta=1}^{\ell+1}F(xv_\beta^{-1}),
\qquad
(S_\ell F)(x)=F(xz_\ell^{-1}),
\]
where $z_\ell$ has component
$\pmat{\ell}{0}{0}{\ell}$ at $\ell$ and identity components
elsewhere.
Thus $S_\ell$ is the operator for the central double coset
\[
\GL_2(\Z_\ell)
\pmat{\ell}{0}{0}{\ell}
\GL_2(\Z_\ell).
\]
These operators commute with $G$ and generate $\T$.
\end{example}

\begin{example}[Modular curves]\label{ex:modular}
Assume $p\ge5$.
Let $K^p\subseteq\GL_2(\A_f^{(p)})$ be compact open, and let
$Y(K^pL)$ be the open modular curve of level $K^pL$.
Homology is Betti homology of the complex curve.
For cohomology we identify Betti cohomology with geometric
\'etale cohomology over $\overline\Q$ by comparison, and use the
latter to define the Galois action.

We fix the level $K^p$ away from $p$ and vary the level at $p$.
For $n\ge1$, put
\[
Y_n=Y(K^pK_n).
\]

The inclusions $K_{n+1}\subseteq K_n$ give covering maps
$Y_{n+1}\to Y_n$.
Pushforward on homology and pullback on cohomology give
\[
H_1(Y_{n+1},\cO)\longrightarrow H_1(Y_n,\cO),
\qquad
H^1(Y_n,\F)\longrightarrow H^1(Y_{n+1},\F).
\]
Thus homology forms an inverse system, while cohomology forms
a direct system.

Let $\m$ be non-Eisenstein, meaning that $\bar\rho_\m$ is
absolutely irreducible.
The completed homology is
\[
\Ht=\varprojlim_n H_1(Y_n,\cO)_\m.
\]
The stabilizers at level $K^p\Iwone$ are finite and inject into
$\Iwone$. Since $p\ge5$, the pro-$p$ group $\Iwone$ is
torsion-free. An element of order $p$ would have a primitive
$p$th root of unity as an eigenvalue, of degree $p-1>2$ over
$\Q_p$. Thus $K^p\Iwone$ is neat.

By homology--cohomology duality and
\cite[Corollary~5.3.19]{Em11}, $\Ht$ is finite free over
$\cO[[\Iwone]]$. Since $[K:\Iwone]=(p-1)^2(p+1)$ is prime
to $p$, averaging over coset representatives shows that $\Ht$
is finitely generated and projective over $\cO[[K]]$.
The tower supplies a continuous right $G$-action on $\Ht$
extending its $K$-action, and the tame Hecke operators commute
with this action.
Thus $\Ht$ is a projective augmented module.
Its residual smooth dual is
\[
\Pi=(\Ht/\varpi)^\vee
\simeq\varinjlim_n H^1(Y_n,\F)_\m.
\]
We use the adelic, Galois, and Hecke actions of
\cite[\S\S5.1--5.2]{Em11}.

Put $J_r^1=J_r\cap\Iwone$ and define
\[
H^1(K^pJ_r;\varepsilon,\F)
=H^1(Y(K^pJ_r^1),\F)^{(\varepsilon)},
\]
where the superscript is the character space for
$J_r/J_r^1=T(\F_p)$.

Since $\m$ is non-Eisenstein, we have
$H^0(Y_n,\F)_\m=0,$
for every $n\ge1$.

Applying Hochschild--Serre to the covers
$Y_n\to Y(K^pJ_r^1)$ for sufficiently large $n$,
and passing to the direct limit, gives
\[
H^1(Y(K^pJ_r^1),\F)_\m
\simeq\Pi^{J_r^1}.
\]
Taking the $\varepsilon$-character spaces for
$J_r/J_r^1=T(\F_p)$ and applying Lemma \ref{lem:level}, we obtain
\[
V_r\simeq H^1(K^pJ_r;\varepsilon,\F)_\m.
\]

For $\varepsilon=1$, this is cohomology at level $K^pJ_r$;
otherwise it is the specified character space on the cover.
Equation \eqref{eq:Jr-conj} and Lemma \ref{lem:level} give the
interpretation as $\Gamma_0(p^{r+1})$-type at $p$.
\end{example}

\section{Square-root nilpotence}
\label{sec:estimates}
In this section, we prove Theorem~\ref{thm:A}, which bounds the
Hecke nilpotence exponent on $V_{r,\m}$ by $O(p^{r/2})$.

For $n\ge1$, recall the principal congruence subgroup
\[
K_n=I_2+p^nM_2(\Z_p)
=\pmat{1+p^n\Z_p}{p^n\Z_p}{p^n\Z_p}{1+p^n\Z_p}.
\]

Fix an occurring maximal ideal $\m\subseteq\T$.
By Lemma \ref{lem:Pi}(1), the invariant space $\Pi^{K_1}$
is finite-dimensional.
Its local component $(\Pi^{K_1})_\m$ is nonzero, since
$V_{1,\m}\ne0$ by Proposition \ref{prop:density} and
$K_1\subseteq J_1$.
We may therefore choose an integer
$1\le b_\m\le\dim_\F(\Pi^{K_1})_\m$ such that
\begin{equation}\label{eq:b}
\m^{b_\m}(\Pi^{K_1})_\m=0.
\end{equation}
Indeed, the nonzero terms of its $\m$-adic filtration
strictly decrease by Nakayama's lemma.

To extend \eqref{eq:b} to vectors fixed by $K_n$, put
$Q_n=K_1/K_n$ and let
\[
I_n=\ker\bigl(\F[Q_n]\longrightarrow\F\bigr)
\]
be its augmentation ideal, where the map sends every group
element to $1$.
The quotient $Q_n$ acts on $(\Pi^{K_n})_\m$, with
\[
\bigl((\Pi^{K_n})_\m\bigr)^{Q_n}
=(\Pi^{K_1})_\m.
\]
Thus $I_nv=0$ means that $v$ is fixed by $K_1$.

We first prove that, for $v\in(\Pi^{K_n})_\m$ and $t\ge1$,
the relation $I_n^tv=0$ implies $\m^{b_\m t}v=0$.

\begin{lemma}\label{lem:transfer}
For every $n,t\ge1$ and $v\in(\Pi^{K_n})_\m$,
\[
I_n^tv=0
\quad\Longrightarrow\quad
\m^{b_\m t}v=0.
\]
\end{lemma}

\begin{proof}
For $j\ge0$, put
\[
M_j=\{v\in(\Pi^{K_n})_\m:I_n^jv=0\}.
\]
Then $M_0=0$ and $M_1=(\Pi^{K_1})_\m$.
For $j\ge1$ and $v\in M_j$, one has
$I_n^{j-1}v\subseteq M_1$.
Since $\T$ commutes with $Q_n$, equation \eqref{eq:b} gives
\[
I_n^{j-1}(\m^{b_\m}v)
=\m^{b_\m}(I_n^{j-1}v)=0.
\]
Thus $\m^{b_\m}M_j\subseteq M_{j-1}$.
Iterating $t$ times gives $\m^{b_\m t}M_t=0$.
\end{proof}

Next, we find an exponent $t$ for which $I_n^tv=0$
when $v$ is fixed by one of the subgroups
\[
H_{a,b}=\pmat{1+p\Z_p}{p^a\Z_p}{p^b\Z_p}{1+p\Z_p},
\qquad a,b\ge1.
\]
Note that $K_n\subseteq H_{a,b}$ whenever $n\ge\max(a,b)$.

\begin{lemma}\label{lem:perm}
Let $a,b\ge1$ and $n\ge\max(a,b)$.
Put $H=H_{a,b}/K_n\subseteq Q_n$.
Then
\[
I_n^{\,p^{a-1}+p^{b-1}-1}\F[Q_n/H]=0.
\]
\end{lemma}

\begin{proof}
For $\smat xyvw\in K_1$, we have the decomposition
\[
\pmat{x}{y}{v}{w}
=
\pmat{1}{0}{v/x}{1}
\pmat{x}{0}{0}{w-yv/x}
\pmat{1}{y/x}{0}{1},
\]
where $y/x,v/x\in p\Z_p$ and $x,w-yv/x\in1+p\Z_p$.
Since $1+p$ topologically generates $1+p\Z_p$ for odd $p$,
the group $K_1$ is topologically generated by
\[
g_1=\pmat{1+p}{0}{0}{1},\quad
g_2=\pmat{1}{0}{0}{1+p},\quad
g_3=\pmat{1}{p}{0}{1},\quad
g_4=\pmat{1}{0}{p}{1}.
\]
Put $\bar g_i=g_iK_n$ and $u_i=\bar g_i-1$.
Since the $\bar g_i$ generate $Q_n$, the augmentation ideal is
\[
I_n=(u_1,u_2,u_3,u_4)\subseteq\F[Q_n].
\]

Let $[K_1,K_1]$ denote the commutator subgroup of $K_1$.
We can show that
\[
[K_1,K_1]\subseteq K_1^p=\{g^p:g\in K_1\}.
\]
Hence, for each $i,j$, we can write
$\bar g_i\bar g_j\bar g_i^{-1}\bar g_j^{-1}=q^p,$ for some $q\in Q_n$.

Since $q-1\in I_n$, characteristic $p$ gives
\[
u_iu_j-u_ju_i
=\bar g_i\bar g_j-\bar g_j\bar g_i
=(q^p-1)\bar g_j\bar g_i
=(q-1)^p\bar g_j\bar g_i
\in I_n^p\subseteq I_n^3.
\]
Thus the algebra
\[
\gr_{I_n}\F[Q_n]
=\bigoplus_{d\ge0}I_n^d/I_n^{d+1}
\]
is generated by the commuting classes $U_i=u_i+I_n^2$.

Let $e_H$ be the basis vector of $\F[Q_n/H]$ corresponding
to the coset $H$.
It generates $\F[Q_n/H]$ as an $\F[Q_n]$-module.
Since
$\bar g_1,\bar g_2,\bar g_3^{p^{a-1}},\bar g_4^{p^{b-1}}\in H$,
one has
\[
u_1e_H=u_2e_H=0,\qquad
u_3^{p^{a-1}}e_H=u_4^{p^{b-1}}e_H=0.
\]
Consequently, $\gr_{I_n}\F[Q_n/H]$ is generated by
\[
\bar e_H=e_H+I_n\F[Q_n/H],
\]
with relations
\[
U_1\bar e_H=U_2\bar e_H=0,\qquad
U_3^{p^{a-1}}\bar e_H=U_4^{p^{b-1}}\bar e_H=0.
\]
Because the $U_i$ commute, this graded module is spanned by
\[
U_3^iU_4^j\bar e_H,
\qquad
0\le i<p^{a-1},\quad 0\le j<p^{b-1}.
\]
These vectors have degree at most $p^{a-1}+p^{b-1}-2$.
Hence, for $t=p^{a-1}+p^{b-1}-1$,
\[
I_n^t\F[Q_n/H]=I_n^{t+1}\F[Q_n/H].
\]
Iterating this equality and using nilpotence of $I_n$
gives $I_n^t\F[Q_n/H]=0$.
\end{proof}

Finally, for $v\in V_{r,\m}$, we apply Lemma \ref{lem:perm}
to $\pi(g_r)v$ for a suitable diagonal matrix $g_r$.
We choose $g_r$ so that the exponent
$t=p^{a-1}+p^{b-1}-1$ is at most $2p^{r/2}$.
Since $\T$ commutes with $G$, a power of $\m$ kills
$\pi(g_r)v$ if and only if it kills $v$.
Lemma \ref{lem:transfer} then gives the desired bound.

\begin{theorem}[Theorem~\ref{thm:A}]\label{thm:balancing}
For every $r\ge0$,
\begin{equation}\label{eq:balancing}
\m^{\,\lceil2b_\m p^{r/2}\rceil}V_{r,\m}=0.
\end{equation}
\end{theorem}

\begin{proof}
For $r=0$, the inclusion $K_1\subseteq J_0=\Iw$ gives
$V_{0,\m}\subseteq(\Pi^{K_1})_\m$, so \eqref{eq:b} proves
\eqref{eq:balancing}.

Now suppose $r\ge1$, put $s=\lfloor r/2\rfloor$ and
$g_r=\pmat{p^{-s}}{0}{0}{1}$.
Conjugation by $g_r$ gives:
\begin{equation}\label{eq:gJg}
g_rJ_rg_r^{-1}
=\pmat{\Z_p^\times}{p^{r-s}\Z_p}{p^{s+1}\Z_p}{\Z_p^\times},
\qquad
g_rJ_rg_r^{-1}\cap K_1=H_{r-s,s+1}.
\end{equation}
Since $1\le r-s\le s+1$, the subgroup $g_rJ_rg_r^{-1}$
contains $K_{s+1}$.

Let $v\in V_{r,\m}$ and put $v'=\pi(g_r)v$.
Diagonal conjugation preserves the diagonal entries, so the
conjugated character $\varepsilon$ is trivial on $H_{r-s,s+1}$.
Since $\T$ commutes with $\pi(g_r)$, the vector $v'$ belongs to
$(\Pi^{K_{s+1}})_\m$ and is fixed by $H_{r-s,s+1}$.

Put $H=H_{r-s,s+1}/K_{s+1}\subseteq Q_{s+1}$.
The map
\[
\F[Q_{s+1}/H]\longrightarrow(\Pi^{K_{s+1}})_\m,
\qquad gH\longmapsto gv',
\]
is $\F[Q_{s+1}]$-linear.
Lemma \ref{lem:perm} therefore gives $I_{s+1}^tv'=0$, where
\[
t=p^{r-s-1}+p^s-1\le2p^{r/2}.
\]
The choice $s=\lfloor r/2\rfloor$ makes both exponents
$r-s-1$ and $s$ at most $r/2$.
Lemma \ref{lem:transfer} gives $\m^{b_\m t}v'=0$.
Applying $\pi(g_r^{-1})$, which commutes with $\T$, gives
$\m^{b_\m t}v=0$.
The inequality $b_\m t\le\lceil2b_\m p^{r/2}\rceil$
proves \eqref{eq:balancing}.
\end{proof}

The images of $\m$ and $\n_\m$ in
$\End_\F(V_{r,\m})$ coincide, so \eqref{eq:balancing}
is the bound asserted in Theorem \ref{thm:A}.

Using $V_r=\Sbar_{\deg\le p^r-1}$ from Section \ref{sec:level},
we obtain the corresponding square-root bound in Mahler degree $N$.

\begin{corollary}\label{cor:nilpotence}
With $C_\m=2b_\m\sqrt p$, one has, for every $N\ge0$,
\begin{equation}\label{eq:A2}
\m^{\,\lceil C_\m\sqrt{N+1}\rceil}\Sbar_{\deg\le N,\m}=0.
\end{equation}
\end{corollary}

\begin{proof}
Choose $r\ge0$ minimal with $N<p^r$.
The identity $V_r=\Sbar_{\deg\le p^r-1}$ established in
Section \ref{sec:level} gives
$\Sbar_{\deg\le N,\m}\subseteq V_{r,\m}$.
Since $p^r\le p(N+1)$, one has
\[
2b_\m p^{r/2}\le C_\m\sqrt{N+1}.
\]
Theorem \ref{thm:balancing} now gives \eqref{eq:A2}.
\end{proof}

\begin{remark}\label{rem:sharper-nilpotence}
The proof of Theorem \ref{thm:balancing} also gives a sharper bound
$\m^{b_\m\lambda_p(r)}V_{r,\m}=0$, where $\lambda_p(0)=1$
and, for $r\ge1$,
\[
\lambda_p(r)=p^{\lfloor(r-1)/2\rfloor}
+p^{\lceil(r-1)/2\rceil}-1.
\]
\end{remark}

\section{Dimension and structure of the completed Hecke algebra}
\label{sec:completed}\
The main goal is Theorem \ref{thm:B}. For $p\ge5$, finite length
of $\Pi[\m]$ implies that $A_\m$ is noetherian and has Krull
dimension at least two. If, in addition,
$\dim_\F\n_\m/\n_\m^2\le2$, then $A_\m\simeq\F[[X,Y]]$.
We prove these assertions in order: block theory gives
noetherianness; the socle inequality and the estimates of
Section \ref{sec:estimates} give the dimension bound; and the
cotangent bound gives the power series structure.
We begin by constructing $A_\m$ and transferring the finite-step
Hecke estimates to it.

\subsection{The completed local Hecke algebra}
\label{sec:completed-def}

Fix an occurring maximal ideal $\m$, and recall the completed
local Hecke algebra $A_\m$ from \eqref{eq:Adef}.
We describe $A_\m$ using the degree filtration and compare
the actions of $\m\subseteq\T$ and the maximal ideal
$\n_\m\subseteq A_\m$ at each finite step.

For $N\ge0$, put
\[
A_{N,\m}
=\im\bigl(\T\to\End_\F(\Sbar_{\deg\le N,\m})\bigr).
\]
Since $V_{r,\m}=\Sbar_{\deg\le p^r-1,\m}$ and the indices
$p^r-1$ are cofinal, we have
\[
A_\m=\varprojlim_N A_{N,\m}.
\]
The transition maps, given by restriction, are surjective.
After omitting any initial zero terms, each $A_{N,\m}$
is a finite local $\F$-algebra with residue field $\F$.
Their inverse limit is local as well: an element with nonzero
residue has compatible inverses at every finite step.
Thus $A_\m$ is a compact local $\F$-algebra.

We give $A_\m$ the inverse-limit topology.
Once noetherianness is established, Lemma \ref{lem:cotangent}
identifies it with the $\n_\m$-adic topology.

For every nonzero finite step, the maps $\T\to A_{N,\m}$
and $A_\m\to A_{N,\m}$ send $\m$ and $\n_\m$, respectively,
onto the maximal ideal of $A_{N,\m}$.
Consequently, for every $N\ge0$ and $q\ge1$,
\begin{equation}\label{eq:same-ideal}
\m^q\Sbar_{\deg\le N,\m}
=\n_\m^q\Sbar_{\deg\le N,\m}.
\end{equation}
In particular, the estimates of Corollary \ref{cor:nilpotence}
also hold with $\n_\m$ in place of $\m$, and the level socle is
\[
V_{r,\m}[\n_\m]=V_{r,\m}[\m]
=\Pi[\m]^{(J_r,\varepsilon)}.
\]

\subsection{$A_\m$ is noetherian}
\label{sec:noeth}
We prove the assertion in part (1) of
Theorem \ref{thm:B} that $A_\m$ is noetherian.

Write $\n_\m$ for the maximal ideal of $A_\m$.
Our goal is to find finitely many elements
$a_1,\ldots,a_s\in\n_\m$ such that $A_\m$ is finite over
the closed $\F$-subalgebra $\F[[a_1,\ldots,a_s]]$.
Since this subalgebra is noetherian, this implies that $A_\m$
is also noetherian.

We want to apply $G$-representation theory while retaining
the faithful action of $A_\m$ on $\Sbar_\m$.
Since $\Sbar_\m$ need not be stable under $G$, we consider
\[
W=\F[G]\Sbar_\m
=\left\{\sum_i g_iv_i:g_i\in G,\ v_i\in\Sbar_\m\right\},
\]
where all sums are finite.
This is a $G$-subrepresentation of $\Pi$.

We extend the $A_\m$-action to $W$ by
\[
a\Bigl(\sum_i g_iv_i\Bigr)=\sum_i g_i(av_i),
\qquad a\in A_\m.
\]
To check that this is well defined, choose $L\ge0$ such that
$v_i\in\Sbar_{\deg\le L,\m}$ for every $i$.
Choose $t\in\T$ whose image in $A_{L,\m}$ agrees with
the image of $a$.
Since the $\T$-action on $\Pi$ commutes with the $G$-action,
\[
\sum_i g_iv_i=0
\quad\Longrightarrow\quad
\sum_i g_i(av_i)
=t\Bigl(\sum_i g_iv_i\Bigr)=0.
\]
The extended $A_\m$-action commutes with $G$ and remains
faithful because its restriction to $\Sbar_\m$ is faithful.
The induced $\T$-action agrees with the
$\T$-action on $\Pi$ restricted to $W$.

We next pass to the compact dual
\[
N=W^\vee.
\]
The $\Iw$-action on $W$ induces a continuous right
$\Iw$-action on $N$.
Thus $N$ is a module over the algebra
$\Lambda=\F[[\Iw]].$

Taking the dual turns the inclusion $W\subseteq\Pi$
into a $\Lambda$-linear quotient
\[
\Ht/\varpi\twoheadrightarrow N.
\]
Since $\Ht/\varpi$ is finitely generated over $\Lambda$,
so is $N$.
The induced $A_\m$-action on $N$ is faithful and commutes
with the $\Lambda$-action.

This commutation makes $\n_\m N$ a $\Lambda$-submodule
of $N$.
Since $\Lambda$ is noetherian, $\n_\m N$ is also finitely
generated over $\Lambda$.
We will use this finiteness to choose finitely many
elements of $\n_\m$.

\begin{lemma}\label{lem:finite-operators}
The $\Lambda$-module $\n_\m N$ is finitely generated.
There exist $a_1,\ldots,a_s\in\n_\m$ such that
\begin{equation}\label{eq:finite-hecke-operators}
\n_\m N=a_1N+\cdots+a_sN.
\end{equation}
In particular, $\n_\m N$ is closed in $N$.
\end{lemma}

\begin{proof}
Since $\Lambda$ is noetherian and $N$ is finitely generated
over $\Lambda$, its submodule $\n_\m N$ is finitely generated.

For each $a\in\n_\m$, the image $aN$ is a $\Lambda$-submodule
because $a$ commutes with $\Lambda$.
The sum
\[
\n_\m N=\sum_{a\in\n_\m}aN
\]
is therefore a finite subsum, giving
\eqref{eq:finite-hecke-operators}.
This sum is the image of the continuous map
\[
N^{\oplus s}\longrightarrow N,
\qquad
(n_1,\ldots,n_s)\longmapsto\sum_i a_in_i.
\]
It is compact and hence closed.
\end{proof}

We now use the finite-length hypothesis on $\Pi[\m]$.
Our next aim is to prove that $W$ is \emph{locally finite}:
every vector of $W$ lies in a $G$-subrepresentation of
finite composition length.
Write $\ell_G(X)$ for the composition length of a
finite-length smooth $G$-representation $X$.

\begin{lemma}\label{lem:local-finiteness}
\label{lem:W-properties}
If $\Pi[\m]$ has finite length as a smooth
$G$-representation, then $W$ is locally finite.
\end{lemma}

\begin{proof}
Put $W_0=0$ and $W_j=W[\n_\m^j]$ for $j\ge1$.
These are $G$-subrepresentations because the $A_\m$-action
commutes with $G$.
We show that they exhaust $W$ and have finite length.

Let $w=\sum_i g_iv_i\in W$.
Choose $L$ sufficiently large that
$v_i\in\Sbar_{\deg\le L,\m}$ for every $i$ and
$A_{L,\m}\ne0$.
The ideal $\ker(A_\m\to A_{L,\m})$ annihilates each
$v_i$ and therefore annihilates $w$.
Since the image of $\n_\m$ in the finite local algebra
$A_{L,\m}$ is nilpotent, some power of $\n_\m$ kills $w$.
Thus
\[
W=\bigcup_{j\ge1}W_j.
\]

The natural map $\T\to A_\m$ sends $\m$ into $\n_\m$.
Compatibility with the original $\T$-action gives
\[
W_1=W[\n_\m]\subseteq\Pi[\m].
\]
Hence $W_1$ has finite length.

Fix the operators in Lemma \ref{lem:finite-operators}
and put $I=(a_1,\ldots,a_s)A_\m$.
Starting from $\n_\m N=IN$, commutativity gives inductively
\[
\n_\m^{j+1}N
=\n_\m^j(IN)
=I(\n_\m^jN)
=I^{j+1}N.
\]
Thus $\n_\m^jN=I^jN$ for every $j\ge1$.
Taking annihilators in $W=N^\vee$ gives
\[
W_j=W[\n_\m^j]=W[I^j].
\]
This equality concerns the action on $N$ and does not
require $I=\n_\m$ as ideals of $A_\m$.

Choose generators $c_1,\ldots,c_{d_j}$ of $I^{j-1}$.
Multiplication by these generators maps $W_j$ into
$W_1^{\oplus d_j}$, with kernel $W_{j-1}$.
It induces a $G$-equivariant injection
\[
W_j/W_{j-1}\hookrightarrow W_1^{\oplus d_j},
\qquad
w+W_{j-1}\longmapsto(c_1w,\ldots,c_{d_j}w).
\]
Induction shows that every $W_j$ has finite length,
with all its irreducible subquotients occurring in $W_1$.
Since $W=\bigcup_jW_j$, the representation $W$ is locally finite.
\end{proof}

We next recall the projective-envelope functor used to
convert the finite $G$-length of $W_1$ into a
finite-dimensional quotient.

For a continuous central character
$\zeta:\Q_p^\times\to\cO^\times$, let $\mathfrak C_\zeta$
be the category anti-equivalent, under Pontryagin duality,
to locally finite smooth $\cO$-torsion
$G$-representations with central character $\zeta$.
Morphisms are continuous $\cO$-linear $G$-equivariant maps.

We use compact right $G$-modules; inversion translates
the results of \cite{PT21} to this convention.
On objects annihilated by $\varpi$, Pontryagin duality
agrees with the $\F$-linear duality used above.

Choose pairwise nonisomorphic absolutely irreducible smooth
$\F$-representations $\sigma_1,\ldots,\sigma_t$ with central
character $\zeta\bmod\varpi$.
Their compact duals have projective envelopes in
$\mathfrak C_\zeta$; see \cite[\S4.1]{PT21}.
Put
\[
P_i\twoheadrightarrow\sigma_i^\vee,
\qquad
\mathcal P=\bigoplus_{i=1}^tP_i,
\qquad
m(Q)=\Hom_{\mathfrak C_\zeta}(\mathcal P,Q).
\]

\begin{theorem}\label{thm:PT}
The functor $m$ is exact and faithful on the full
subcategory of compact modules $Q\in\mathfrak C_\zeta$
whose smooth duals have irreducible subquotients among
$\sigma_1,\ldots,\sigma_t$.
It commutes with inverse limits in this subcategory.

For such $Q$, if $Q$ has finite length and $\varpi Q=0$, then
\[
\dim_\F m(Q)=\ell_G(Q^\vee)<\infty.
\]
\end{theorem}

\begin{proof}
Projectivity gives exactness.
If $Q\ne0$, its locally finite smooth dual contains an
irreducible subrepresentation.
Thus $Q$ has a simple quotient
$Q\twoheadrightarrow\sigma_i^\vee$.
Projectivity lifts $P_i\to\sigma_i^\vee$ to a nonzero map
$P_i\to Q$.
Hence $m(Q)\ne0$, and exactness implies faithfulness.

The inverse-limit assertion follows from the universal
property of $\Hom(\mathcal P,-)$.
The projective-envelope property and absolute
irreducibility give
\[
\Hom_{\mathfrak C_\zeta}(P_i,\sigma_j^\vee)
\simeq
\begin{cases}
\F,&i=j,\\
0,&i\ne j.
\end{cases}
\]
Applying exactness to a composition series proves
the dimension formula.
\end{proof}

We will make a finite extension of $\F$ when applying
this construction.
The following lemma identifies the resulting Hecke algebra.

\begin{lemma}\label{lem:basechange}
Let $\F'/\F$ be finite.
For every nonzero $A_{L,\m}$,
\[
A_{L,\m}\otimes_\F\F'
=
\im\bigl(\T\otimes_\F\F'\to
\End_{\F'}(\Sbar_{\deg\le L,\m}\otimes_\F\F')\bigr).
\]
This algebra is local with residue field $\F'$, and
\[
A_\m\otimes_\F\F'
=
\varprojlim_L(A_{L,\m}\otimes_\F\F').
\]
\end{lemma}

\begin{proof}
Scalar extension is exact and commutes with
finite-dimensional endomorphism spaces.
The extended finite algebra has a nilpotent ideal with
quotient $\F'$, so it is local.

Since $\F'/\F$ is finite, tensoring with $\F'$ is a finite
direct sum as a vector-space functor and commutes with
inverse limits.
\end{proof}

We now prove the noetherianness assertion in
Theorem \ref{thm:B}.

\begin{proposition}\label{prop:FLnoeth}
If $\Pi[\m]$ has finite length as a smooth
$G$-representation, then $A_\m$ is noetherian.
\end{proposition}

\begin{proof}
Put $W_j=W[\n_\m^j]$ for $j\ge1$.
By Lemma \ref{lem:finite-operators} and compact duality,
\[
\n_\m N=\sum_{i=1}^s a_iN,
\qquad
(N/\n_\m N)^\vee\simeq W_1.
\]
The representation $W_1\subseteq\Pi[\m]$ has finite
$G$-length, but may be infinite-dimensional over $\F$.
We use the projective-envelope functor of
Theorem \ref{thm:PT} to turn this finite $G$-length into
a finite-dimensional quotient, while retaining a faithful
Hecke action.

By the proof of Lemma \ref{lem:local-finiteness},
$W=\bigcup_jW_j$, each $W_j$ has finite length, and
all irreducible subquotients occur among the finitely many
types $\tau_1,\ldots,\tau_t$ in $W_1$.
Before applying Theorem \ref{thm:PT}, we need two
further properties:
\begin{enumerate}
\item
Each $\tau_i$ is absolutely irreducible with an
$\F^\times$-valued central character.
\item
After modifying the $G$-action, $W$ decomposes into
$A_\m$-stable summands with fixed central characters.
\end{enumerate}

For (1), we extend the coefficient field.
Choose $\F[G]$-generators $v_1,\ldots,v_e$ of $W_1$,
fixed by a common compact open subgroup $U_0$.
Evaluation gives
\[
\End_G(W_1)\hookrightarrow(W_1^{U_0})^e,
\qquad
f\longmapsto(f(v_1),\ldots,f(v_e)).
\]
Since $W_1\subseteq\Pi$ is admissible, the target is
finite-dimensional.
Hence
\[
C=\im\bigl(\F[Z(G)]\to\End_G(W_1)\bigr)
\]
is a finite-dimensional commutative $\F$-algebra.

After a finite extension splitting the residue fields
of $C$, the irreducible constituents have
$\F^\times$-valued central characters.
A further finite extension makes them absolutely
irreducible \cite[\S4.1]{PT21}.
By Lemma \ref{lem:basechange}, we make these extensions,
relabel the $\tau_i$, and retain the same notation.
The equality \eqref{eq:finite-hecke-operators} remains valid.

For (2), we modify the $G$-action while preserving the
Hecke and Iwahori actions.
The compact centre already acts on $W$ through
$zI_2\mapsto\varepsilon(\bar z,\bar z)$, so it suffices
to make $z_0=pI_2$ act semisimply.
Although $z_0$ acts scalarly on each $\tau_i$, it may
act nonsemisimply on $W$.
We remove this unipotent part as follows.

Let $\pi$ denote the $G$-action and put
\[
h=|\F|-1,\qquad U=\pi(z_0)^h.
\]
Every scalar in $\F^\times$ has $h$-th power $1$.
Thus $U$ acts as the identity on every irreducible factor,
and
\[
(U-1)^{\ell_G(W_j)}W_j=0.
\]
If $p^k\ge\ell_G(W_j)$, then $U^{p^k}|_{W_j}=1$.
For $\alpha\in\Z_p$, define $U^\alpha$ on $W_j$ by
reducing $\alpha$ modulo such a $p^k$.
The definition is independent of sufficiently large $k$
and compatible as $j$ varies.

Since $2h\in\Z_p^\times$, define
\[
\pi'(g)=U^{-v_p(\det g)/(2h)}\pi(g).
\]
The operator $U$ commutes with $G$ and $A_\m$, and
$v_p\circ\det$ is additive.
Thus $\pi'$ is a representation commuting with $A_\m$.
It agrees with $\pi$ on compact subgroups and is therefore
smooth.
Since $v_p(\det z_0)=2$, we have
\[
\pi'(z_0)^h=U^{-1}\pi(z_0)^h=1.
\]

On each $W_j$, the operators $U^\alpha$ are polynomials
in $U$.
They preserve an original composition series and act
trivially on its irreducible factors.
Hence the modified action preserves the lengths and
irreducible constituents of the $W_j$.

Write $W^{\mathrm{tw}}=(W,\pi')$.
Since $T^h-1$ splits over $\F$ with distinct roots,
\[
W^{\mathrm{tw}}
=
\bigoplus_{\lambda\in\F^\times}
\ker\bigl(\pi'(z_0)-\lambda\bigr).
\]
These summands are $A_\m$-stable.
On each nonzero summand, the compact centre acts through
$\varepsilon$ and $z_0$ acts through $\lambda$.
Since $\Q_p^\times=p^\Z\times\Z_p^\times$, each summand
has a fixed central character.

Equip $N$ with the dual modified action.
Its $A_\m$-module structure and $\Iw$-action are unchanged,
and
\[
(N/\n_\m N)^\vee
\simeq W^{\mathrm{tw}}[\n_\m].
\]
The representation on the right has the same finite
$G$-length as $W_1$.

With (1) and (2) established, we apply the projective-envelope
functor.
For each central character occurring in $W^{\mathrm{tw}}$,
take its Teichm\"uller lift and the projective envelopes
of the irreducible types with that character.
Apply Theorem \ref{thm:PT} to the corresponding summands
of $N$, and write $m$ for the direct sum of these functors.

The $A_\m$-action on $m(N)$ is given by postcomposition.
Faithfulness of $m$ gives
\[
\Ann_{A_\m}(m(N))=\Ann_{A_\m}(N)=0.
\]
Exactness and the dimension formula give
\[
\begin{aligned}
\dim_\F\frac{m(N)}{\sum_i a_i m(N)}
=\dim_\F m(N/\n_\m N)
=\ell_G(W_1)<\infty.
\end{aligned}
\]

To pass from this finite-dimensional quotient to finite
generation, we also need compactness.
Since $W=\bigcup_jW_j$, Theorem \ref{thm:PT} gives
\[
N=\varprojlim_jW_j^\vee,
\qquad
m(N)=\varprojlim_jm(W_j^\vee),
\]
where the $W_j$ carry the modified action.
The spaces $m(W_j^\vee)$ are finite-dimensional, and
the transition maps are surjective.

The kernel calculation in the proof of
Lemma \ref{lem:local-finiteness} shows that every vector
of $W$ has an open annihilator in $A_\m$.
Each $W_j$ has finitely many $\F[G]$-generators.
Since the $A_\m$-action commutes with $G$, the intersection
of their annihilators kills $W_j$.
Thus $\Ann_{A_\m}(W_j)$ is open, and the action on
$m(W_j^\vee)$ factors through a finite quotient of $A_\m$.
It follows that $m(N)$ is a faithful profinite
$A_\m$-module.

We now prove finite generation over
$\F[[a_1,\ldots,a_s]]$.
The images of the $a_i$ are nilpotent in every finite
quotient of $A_\m$, so substitution gives
\[
\F[[X_1,\ldots,X_s]]
\twoheadrightarrow\F[[a_1,\ldots,a_s]],
\qquad X_i\longmapsto a_i.
\]
This closed subalgebra is therefore compact and noetherian,
with maximal ideal generated by $a_1,\ldots,a_s$.

Choose $\mu_1,\ldots,\mu_d\in m(N)$ lifting an $\F$-basis
of $m(N)/\sum_i a_i m(N)$.
Their images generate
\[
m(W_j^\vee)\big/\sum_i a_i m(W_j^\vee)
\]
over $\F$.
Nakayama's lemma therefore shows that they generate
each $m(W_j^\vee)$ over $\F[[a_1,\ldots,a_s]]$.
Consequently, the continuous map
\[
\F[[a_1,\ldots,a_s]]^d\longrightarrow m(N),
\qquad
(f_1,\ldots,f_d)\longmapsto\sum_i f_i\mu_i
\]
has dense image.
Its image is compact and hence closed, so it is surjective.
Thus $m(N)$ is finite over $\F[[a_1,\ldots,a_s]]$.

Finally, an element of $A_\m$ annihilating all the $\mu_i$
annihilates $m(N)$ and is therefore zero by faithfulness.
We obtain an $\F[[a_1,\ldots,a_s]]$-linear injection
\[
A_\m\hookrightarrow m(N)^d,
\qquad
a\longmapsto(a\mu_1,\ldots,a\mu_d).
\]
The target is finite over the noetherian ring
$\F[[a_1,\ldots,a_s]]$.
Its submodule $A_\m$ is therefore finite over this ring,
and hence noetherian.

Let $\F_0$ be the original coefficient field and
$A_{\m,0}$ the corresponding completed local algebra.
By Lemma \ref{lem:basechange}, the algebra considered
after scalar extension is
\[
A_{\m,0}\otimes_{\F_0}\F.
\]
The map
$A_{\m,0}\longrightarrow A_{\m,0}\otimes_{\F_0}\F$
is faithfully flat.
Noetherianness therefore descends to $A_{\m,0}$.
\end{proof}

\subsection{Socles and the dimension bound}
\label{sec:duality}

Throughout this subsection, assume $p\ge5$.

We prove the dimension bound and the quadratic lower bound
in Theorem \ref{thm:B}.
We compare each finite level space with its socle
\[
V_{r,\m}[\n_\m]=\Pi[\m]^{(J_r,\varepsilon)}.
\]

\begin{lemma}\label{lem:socle}
For $r\ge0$ and $q\ge1$, if $\n_\m^qV_{r,\m}=0$, then
\begin{equation}\label{eq:socle-bound}
\dim_\F V_{r,\m}
\le
\dim_\F V_{r,\m}[\n_\m]\cdot
\dim_\F(A_\m/\n_\m^q).
\end{equation}
\end{lemma}

\begin{proof}
The case $V_{r,\m}=0$ is immediate.
Equip the dual $V_{r,\m}^\vee=\Hom_\F(V_{r,\m},\F)$ with the
$A_\m$-action $(a\ell)(v)=\ell(av)$.
The evaluation pairing gives
\[
\bigl(V_{r,\m}^\vee/\n_\m V_{r,\m}^\vee\bigr)^\vee
\simeq V_{r,\m}[\n_\m].
\]
By Nakayama's lemma, $V_{r,\m}^\vee$ is generated over
$A_\m$ by $\dim_\F V_{r,\m}[\n_\m]$ elements.
Since $\n_\m^q$ also annihilates $V_{r,\m}^\vee$,
there is a surjection
\[
(A_\m/\n_\m^q)^{\oplus\dim_\F V_{r,\m}[\n_\m]}
\twoheadrightarrow V_{r,\m}^\vee.
\]
Taking dimensions proves \eqref{eq:socle-bound}.
\end{proof}
To apply Lemma \ref{lem:socle}, we bound the dimensions
of the level socles using the finite-length hypothesis
on $\Pi[\m]$.

\begin{proposition}\label{prop:level-socle}
Assume $p\ge5$ and that $\Pi[\m]$ has finite length
as a smooth $G$-representation.
There is a constant $c_0>0$ such that
\begin{equation}\label{eq:linear-socle-growth}
\dim_\F\Pi[\m]^{(J_r,\varepsilon)}
\le c_0(r+1)\textrm{ for all }r\ge0.
\end{equation}
\end{proposition}

\begin{proof}
After a finite scalar extension, we may choose
a composition series of $\Pi[\m]$ with absolutely
irreducible factors $\tau_1,\ldots,\tau_t$ having central
characters, as in the proof of Proposition \ref{prop:FLnoeth}.

For $n\ge1$, put
\[
T_1(p^n)
=\pmat{1+p\Z_p}{p^n\Z_p}{p^n\Z_p}{1+p\Z_p}.
\]
Hu's theorem \cite[Theorem 4.1]{Hu21} applies to each
$\tau_j$ with its own central character.
Left exactness of invariants gives
\[
\dim_\F\Pi[\m]^{T_1(p^n)}
\le
\sum_{j=1}^t\dim_\F\tau_j^{T_1(p^n)}
\le cn \textrm{ for all } n\ge1
\]
for some $c>0$.
Thus a fixed central character on the whole of $\Pi[\m]$
is not needed.

For $r\ge1$, put $s=\lfloor r/2\rfloor$ and
$g_r=\pmat{p^{-s}}{0}{0}{1}$.
By \eqref{eq:gJg},
\[
T_1(p^{s+1})\subseteq g_rJ_rg_r^{-1}\cap K_1.
\]
Since $\varepsilon$ is trivial on this subgroup,
conjugation by $g_r$ gives
\[
\dim_\F\Pi[\m]^{(J_r,\varepsilon)}
\le\dim_\F\Pi[\m]^{T_1(p^{s+1})}
\le c(s+1)
\le c(r+1).
\]
The case $r=0$ follows from admissibility after increasing
the constant.
Finite scalar extension commutes with invariants and
preserves their dimensions, so
\eqref{eq:linear-socle-growth} holds over the original
field $\F$.
\end{proof}

We now combine this socle bound with Theorem \ref{thm:A}
and the noetherianness proved in Proposition \ref{prop:FLnoeth}.

\begin{theorem}
\label{thm:criterion}\label{cor:dim2}
Let $p\ge5$.
Suppose that $\Pi[\m]$ has finite length as a smooth
$G$-representation.
Then $\dim A_\m\ge2$, and there is a constant $c_\m>0$ such that
\[
\dim_\F(A_\m/\n_\m^q)\ge c_\m q^2
\]
for all sufficiently large integers $q$.
\end{theorem}

\begin{proof}
By Proposition \ref{prop:FLnoeth}, $A_\m$ is noetherian.
Put $q_r=\lceil2b_\m p^{r/2}\rceil.$
Theorem \ref{thm:A} gives $\n_\m^{q_r}V_{r,\m}=0$.
For sufficiently large $r$, \eqref{eq:level-count},
Lemma \ref{lem:socle}, and
\eqref{eq:linear-socle-growth} give
\[
\frac{\rho_\m}{2}p^r
\le\dim_\F V_{r,\m}
\le c_0(r+1)\dim_\F(A_\m/\n_\m^{q_r}).
\]
Since $q_r\le(2b_\m+1)p^{r/2}$, it follows that
\[
\frac{\dim_\F(A_\m/\n_\m^{q_r})}{q_r}
\ge
\frac{\rho_\m}{2c_0(2b_\m+1)}
\frac{p^{r/2}}{r+1}
\longrightarrow\infty.
\]

The residue field of $A_\m$ is $\F$, so
\[
\dim_\F(A_\m/\n_\m^q)
=
\operatorname{length}_{A_\m}(A_\m/\n_\m^q).
\]
By Hilbert--Samuel theory \cite[Theorem 11.14]{AM69},
this agrees, for sufficiently large $q$, with a polynomial
of degree $\dim A_\m$ and positive leading coefficient.
If $\dim A_\m\le1$, then
$\dim_\F(A_\m/\n_\m^q)=O(q)$, contradicting the limit above.
Hence $\dim A_\m\ge2$.
The same polynomial gives
$\dim_\F(A_\m/\n_\m^q)\ge c_\m q^2$
for some $c_\m>0$ and all sufficiently large integers $q$.
\end{proof}

\begin{remark}
The finite-length hypothesis cannot be omitted in general.
For a scalar Hecke action on a nonzero residual space,
$A_\m=\F$ and $\dim A_\m=0$.
\end{remark}

\subsection{Power series structure}
\label{sec:structure}

Under the finite-length hypothesis, Proposition
\ref{prop:FLnoeth} gives noetherianness of $A_\m$, and
Theorem \ref{thm:criterion} gives $\dim A_\m\ge2$.
To finish Theorem \ref{thm:B}, we construct a power series
presentation from cotangent generators and then use the
dimension bound to show that its kernel is zero.

Let $x_1,\ldots,x_d\in\n_\m$ have images spanning
$\n_\m/\n_\m^2$ over $\F$.
Their images lie in the nilpotent maximal ideal of every
nonzero finite quotient $A_{L,\m}$.
Thus substitution defines compatible maps to these quotients,
giving a continuous map
\[
\Phi:\F[[X_1,\ldots,X_d]]\longrightarrow A_\m,
\qquad X_i\longmapsto x_i.
\]

\begin{lemma}\label{lem:cotangent}
The map $\Phi$ is surjective.
Moreover, the inverse-limit topology on $A_\m$ agrees
with its $\n_\m$-adic topology.
\end{lemma}

\begin{proof}
We first prove surjectivity at each finite level.
The images of the $x_i$ span the maximal ideal of
$A_{L,\m}$ modulo its square.
Nakayama's lemma shows that they generate this maximal ideal.
Its nilpotence then implies that every element of
$A_{L,\m}$ is a polynomial in their images.
Hence $\Phi$ surjects onto every $A_{L,\m}$.

The image of $\Phi$ is therefore dense in $A_\m$.
It is also compact and hence closed.
Thus $\Phi$ is surjective.

Put $\mathfrak q=(X_1,\ldots,X_d)$ and $K=\ker\Phi$.
Since $\Phi$ is surjective and induces the identity
on the residue field,
\[
\Phi(\mathfrak q)=\n_\m.
\]
Consequently,
\[
\Phi(\mathfrak q^q)=\n_\m^q,
\qquad
\Phi^{-1}(\n_\m^q)=\mathfrak q^q+K
\qquad(q\ge1).
\]
The map $\Phi$ is a quotient map, since its source is compact
and its target is Hausdorff.
Thus each $\n_\m^q$ is open.
Conversely, if $U$ is an open neighbourhood of zero in $A_\m$,
then $\Phi^{-1}(U)$ contains $\mathfrak q^q$ for some $q$,
so $U$ contains $\n_\m^q$.
Hence the inverse-limit topology agrees with the
$\n_\m$-adic topology.
\end{proof}

\begin{corollary}\label{cor:power-series}
Let $p\ge5$, and suppose that $\Pi[\m]$ has finite length
as a smooth $G$-representation.
If $\dim_\F(\n_\m/\n_\m^2)\le2$, then $A_\m\simeq\F[[X,Y]]$.
\end{corollary}

\begin{proof}
Choose $x,y\in\n_\m$ whose images span
$\n_\m/\n_\m^2$.
Applying Lemma \ref{lem:cotangent} with these generators gives
a continuous surjection
\[
\Phi:\F[[X,Y]]\twoheadrightarrow A_\m,
\qquad X\longmapsto x,\quad Y\longmapsto y.
\]
It remains to show that $\ker\Phi=0$.

By Theorem \ref{thm:criterion}, $\dim A_\m\ge2$.
Since $\F[[X,Y]]$ is a two-dimensional domain, a quotient
by a nonzero ideal has dimension at most one.
Thus $\ker\Phi\ne0$ would imply
\[
\dim A_\m
=
\dim\bigl(\F[[X,Y]]/\ker\Phi\bigr)
\le1,
\]
a contradiction.
Hence $\ker\Phi=0$, and $\Phi$ is an isomorphism of
complete local $\F$-algebras.
\end{proof}

\section{Applications to modular curves and definite quaternion algebras}
\label{sec:arith}

Throughout this section, assume $p\ge5$.
We apply Theorem \ref{thm:B} to the two arithmetic realizations of
$\Pi$ from Section \ref{sec:examples}.
In each case, finite length of the Hecke eigenspace $\Pi[\m]$
gives noetherianness and the lower bound $\dim A_\m\ge2$.
Under additional hypotheses, a Galois deformation quotient
identifies $A_\m$ with a power series ring in two variables.

In either setting, fix a non-Eisenstein maximal ideal $\m$
and a character $\varepsilon$ with
$\Sbar^{(\varepsilon)}_\m\ne0$.
Write
\[
\bar\rho=\bar\rho_\m:
G_{\Q,\mathcal R}\longrightarrow\GL_2(\F)
\]
for the associated absolutely irreducible odd residual representation,
where $\mathcal R$ is specified in each case.
With arithmetic Frobenius, our normalization is
\[
\operatorname{tr}\bar\rho(\Frob_\ell)=T_\ell\bmod\m,
\qquad
\det\bar\rho(\Frob_\ell)=\ell S_\ell\bmod\m
\qquad(\ell\notin\mathcal R).
\]
\subsection{Definite quaternion algebras}
\label{sec:definite}

Let $D/\Q$ be definite and split at $p$, with tame level $K^p$
as in Example~\ref{ex:quaternion}. Let $\mathcal R$ contain $p$,
the primes where $D$ ramifies, and the primes where $K^p$ is
not maximal. Put
\[
X=D^\times\backslash D_f^\times/K^p,
\qquad
\Ht=\cO[[X]],
\qquad
\Pi=C^\infty(X,\F).
\]
The splitting at $p$ gives the action of $G=\GL_2(\Q_p)$ on
$\Pi$. We assume that $K^p\Iw$ is neat, so $\Ht$ is projective
over $\cO[[\Iw]]$ by Example~\ref{ex:quaternion}.

Fix a non-Eisenstein eigensystem $\m$ occurring in
$\Sbar^{(\varepsilon)}$. If $\Pi[\m]$ has finite length as a
smooth $G$-representation, Theorem~\ref{thm:B} gives
$\dim A_\m\ge2$. To identify $A_\m$ with $\F[[X,Y]]$, we will
construct a surjection from $\F[[X,Y]]$ onto $A_\m$ under
additional arithmetic hypotheses.

Let
\[
\bar\rho=\bar\rho_\m\colon
G_{\Q,\mathcal R}\longrightarrow\GL_2(\F)
\]
be the absolutely irreducible odd residual Galois representation
attached to $\m$. For $\ell\notin\mathcal R$, it satisfies
\[
\operatorname{tr}\bar\rho(\Frob_\ell)=T_\ell\bmod\m,
\qquad
\det\bar\rho(\Frob_\ell)=\ell S_\ell\bmod\m.
\]
The map onto $A_\m$ will come from a lift of $\bar\rho$ over
$A_\m$. We construct this lift from Galois representations over
the integral Hecke algebras at finite level.

For $r\ge0$, put
\[
M_r=\Map(X/J_r^1,\cO),
\qquad
\mathbb T^D_r
=\cO[T_\ell,S_\ell\mid\ell\notin\mathcal R]
\subseteq\End_\cO(M_r).
\]
The eigensystem $\m$ occurs in
$M_r/\varpi M_r\simeq\Pi^{J_r^1}$ because $\Pi[\m]\ne0$ and
$J_r^1$ is pro-$p$. Let $\widetilde{\m}_r$ be the kernel of
the residual eigensystem
\[
\mathbb T^D_r\longrightarrow\F,
\qquad
T_\ell\longmapsto T_\ell\bmod\m,
\qquad
S_\ell\longmapsto S_\ell\bmod\m,
\]
and put $\mathbb T^D_{r,\m}=(\mathbb T^D_r)_{\widetilde{\m}_r}$.

\begin{lemma}\label{lem:JL}
For each $r\ge0$, there is a continuous lift
\[
\rho_r\colon G_{\Q,\mathcal R}\longrightarrow
\GL_2(\mathbb T^D_{r,\m})
\]
of $\bar\rho$ such that, for $\ell\notin\mathcal R$,
\[
\operatorname{tr}\rho_r(\Frob_\ell)=T_\ell,
\qquad
\det\rho_r(\Frob_\ell)=\ell S_\ell.
\]
\end{lemma}

\begin{proof}
Write $T_\ell^+,S_\ell^+$ for the forward double-coset operators.
For $K_\ell=\GL_2(\Z_\ell)$ and $a_\ell=\smat\ell001$, we have
\[
K_\ell a_\ell^{-1}K_\ell=z_\ell^{-1}K_\ell a_\ell K_\ell.
\]
Thus the operators of Example~\ref{ex:quaternion} satisfy
\[
T_\ell=T_\ell^+(S_\ell^+)^{-1},
\qquad
S_\ell=(S_\ell^+)^{-1}.
\]
We use the same formulas for the classical operators
$T_\ell,S_\ell$. Since the $S_\ell^+$ have finite order in
weight two, the two sets of operators generate the same
integral Hecke algebra.

By Jacquet--Langlands \cite{JL70}, every characteristic-zero
eigensystem of the non-Eisenstein Hecke algebra
$\mathbb T^D_{r,\m}$ occurs on a classical weight-two cusp form
with the same $T_\ell,S_\ell$ eigenvalues for
$\ell\notin\mathcal R$. The finitely many transfers occur at
a common level maximal outside $\mathcal R$. Let
$\mathbb T^{\mathrm{cl}}_{r,\m}$ be the integral tame Hecke
algebra at that level localized at $\m$.

Every relation among the classical $T_\ell,S_\ell$ therefore
vanishes on all quaternionic Hecke eigenspaces. By semisimplicity
in characteristic zero and torsion-freeness over $\cO$, the
relation holds in $\mathbb T^D_{r,\m}$. Sending $T_\ell,S_\ell$
to the corresponding operators gives a surjection
\[
\mathbb T^{\mathrm{cl}}_{r,\m}\twoheadrightarrow
\mathbb T^D_{r,\m}.
\]
The forward residual representation is $\bar\rho^\vee\otimes\omega$
by the Frobenius identities, so it is absolutely irreducible.
Let $\rho^+$ be the classical Galois representation over
$\mathbb T^{\mathrm{cl}}_{r,\m}$ with arithmetic Frobenius
polynomial $X^2-T_\ell^+X+\ell S_\ell^+$ \cite{Ca94}.
At $\Frob_\ell$, the representation
$(\rho^+)^\vee\otimes\chi_{\mathrm{cyc}}$ has trace
$T_\ell^+(S_\ell^+)^{-1}=T_\ell$ and determinant
$\ell(S_\ell^+)^{-1}=\ell S_\ell$.
Base change along this surjection gives $\rho_r$.
Its Frobenius identities identify its reduction with $\bar\rho$.
\end{proof}

The finite quotients of $A_\m=\varprojlim_N A_{N,\m}$ come
from the Hecke action on the Mahler steps
$\Sbar_{\deg\le N,\m}$. Each such step lies in the reduction
of $M_r$ after localization for sufficiently large $r$.
The lift $\rho_r$ can therefore be specialized to $A_{N,\m}$.

\begin{proposition}\label{prop:noeth}
There is a continuous lift
\[
\rho\colon G_{\Q,\mathcal R}\longrightarrow\GL_2(A_\m)
\]
of $\bar\rho$ such that, for $\ell\notin\mathcal R$,
\[
\operatorname{tr}\rho(\Frob_\ell)=T_\ell,
\qquad
\det\rho(\Frob_\ell)=\ell S_\ell.
\]
\end{proposition}

\begin{proof}
For each nonzero Mahler step, choose $r$ with $N<p^r$.
The inclusions
\[
\Sbar_{\deg\le N,\m}\subseteq V_{r,\m}
\subseteq(\Pi^{J_r^1})_\m
\simeq(M_r/\varpi M_r)_\m
\]
give a surjection
\[
\mathbb T^D_{r,\m}/\varpi\twoheadrightarrow A_{N,\m}.
\]
Base change of $\rho_r$ gives a lift over $A_{N,\m}$.

The lifts agree as deformation classes under the transition
maps. Indeed, after reduction to a common quotient, their
traces at $\Frob_\ell$ are the image of the same operator
$T_\ell$. Chebotarev gives equality of traces everywhere,
and absolute irreducibility of $\bar\rho$ then gives equality
of the deformation classes \cite{Ca94}. The resulting maps
\[
R_{\bar\rho,\mathcal R}/\varpi\longrightarrow A_{N,\m}
\]
from the universal deformation ring \cite{Ma89} are compatible.
Their inverse limit gives the map to $A_\m$ that defines $\rho$.
\end{proof}

A fixed-determinant deformation ring will give the source
with two variables. To obtain a map from this ring, we must
fix the determinant of $\rho$. The formula
$\det\rho(\Frob_\ell)=\ell S_\ell$ reduces this to controlling
the action of central scalars.

Embed the finite ideles as scalars in $D_f^\times$ and put
\[
K_{Z^p}=\{z\in(\A_f^{(p)})^\times\mid z\in K^p\},
\qquad
U_{p,\varepsilon}
=\{z\in\Z_p^\times\mid\varepsilon(\bar z,\bar z)=1\}.
\]
Rational scalars and scalars in $K^p$ act trivially, and the
compact centre at $p$ acts through $\varepsilon$. Thus the
scalar action on $\Sbar^{(\varepsilon)}_\m$ factors through
\[
C_\varepsilon
=\A_f^\times/(\Q^\times K_{Z^p}U_{p,\varepsilon}).
\]
This group is finite because
$\A_f^\times=\Q^\times\widehat\Z^\times$ and
$K_{Z^p}U_{p,\varepsilon}$ is open in $\widehat\Z^\times$.
Let $\Gamma$ be the Sylow $p$-subgroup of $C_\varepsilon$.

Let $\chi_{\mathrm{cyc}}$ be the $p$-adic cyclotomic character
and $\omega$ its reduction. Put
\[
\eta=\chi_{\mathrm{cyc}}\,
\mathrm{Teich}(\det\bar\rho\,\omega^{-1}),
\]
and let $R^\eta_{\bar\rho,\mathcal R}$ be the universal
deformation ring with determinant $\eta$ and no ramification
outside $\mathcal R$.

\begin{proposition}\label{prop:upper}
If $\Gamma=1$, there is a continuous surjection
\[
R^\eta_{\bar\rho,\mathcal R}/\varpi\twoheadrightarrow A_\m.
\]
\end{proposition}

\begin{proof}
By class field theory, the scalar action gives a character
$\psi_Z$ of $G_{\Q,\mathcal R}$ through $C_\varepsilon$ with
$\psi_Z(\Frob_\ell)=S_\ell$ for $\ell\notin\mathcal R$.
By Chebotarev, every value of $\psi_Z$ occurs as one of these
$S_\ell$, so $\psi_Z$ takes values in $A_\m^\times$.
Comparing determinants at Frobenius elements gives
\[
\det\rho=\omega\psi_Z.
\]

The ratio $\psi_Z/\overline{\psi_Z}$ takes values in the
pro-$p$ group $1+\n_\m$, where we view $\overline{\psi_Z}$
as $A_\m^\times$-valued via $\F\subset A_\m$. Since
$\Gamma=1$, the ratio also has finite order prime to $p$.
It is therefore trivial, and
\[
\det\rho=\omega\overline{\psi_Z}=\eta\bmod\varpi.
\]
The universal property now gives a continuous map
$R^\eta_{\bar\rho,\mathcal R}/\varpi\to A_\m$.
Its image contains the $T_\ell$ as Frobenius traces.
It also contains the $S_\ell$, which are scalars in $\F$
since $\psi_Z=\overline{\psi_Z}$. The image is closed by
compactness, and these operators topologically generate
$A_\m$, so the map is surjective.
\end{proof}

We call $\bar\rho$ \emph{unobstructed relative to $\mathcal R$}
if
\[
H^2(G_{\Q,\mathcal R},\ad^0\bar\rho)=0,
\]
where $\ad^0\bar\rho$ is the space of trace-zero endomorphisms
with the conjugation action.

\begin{corollary}\label{cor:arithmetic}
Suppose that $\Pi[\m]$ has finite length as a smooth
$G$-representation, that $\Gamma=1$, and that $\bar\rho$
is unobstructed relative to $\mathcal R$. Then
$A_\m\simeq\F[[X,Y]]$.
\end{corollary}

\begin{proof}
Unobstructedness makes $R^\eta_{\bar\rho,\mathcal R}$ formally
smooth over $\cO$. Its tangent space
$H^1(G_{\Q,\mathcal R},\ad^0\bar\rho)$ has dimension two by
the global Euler characteristic formula, using oddness and
$H^0=H^2=0$. Hence
\[
R^\eta_{\bar\rho,\mathcal R}/\varpi\simeq\F[[X,Y]].
\]
Proposition~\ref{prop:upper} gives a surjection from this ring
onto $A_\m$. Its kernel is zero because a nonzero quotient
ideal in $\F[[X,Y]]$ would give $\dim A_\m\le1$, whereas
Theorem~\ref{thm:B} gives $\dim A_\m\ge2$.
\end{proof}
\subsection{Modular curves}
\label{sec:modular-curves}

We use the modular-curve tower of Example~\ref{ex:modular}.
Let $\Sigma_0$ be the set of primes where $K^p$ is not maximal,
and put $\mathcal R=\Sigma_0\cup\{p\}$. Write $K_{\Sigma_0}$
for the factor of $K^p$ at these primes. With $Y_n=Y(K^pK_n)$,
the identifications in that example give
\[
\Pi\simeq\varinjlim_n H^1(Y_n,\F)_\m,
\qquad
V_{r,\m}\simeq H^1(K^pJ_r;\varepsilon,\F)_\m.
\]
To apply Theorem~\ref{thm:B}, we use the restriction
$\bar\rho_p=\bar\rho|_{G_{\Q_p}}$ of the residual representation
attached to $\m$ to establish finite length of $\Pi[\m]$.

\begin{proposition}\label{prop:finite-length}
Suppose that, over $\overline\F$, the representation $\bar\rho_p$
is isomorphic to neither
\[
\chi\otimes\begin{pmatrix}1&*\\0&1\end{pmatrix}
\qquad\text{nor}\qquad
\chi\otimes\begin{pmatrix}1&*\\0&\omega\end{pmatrix}
\]
for any continuous character $\chi:G_{\Q_p}\to\overline\F^\times$,
with $*$ allowed to be zero. Then $\Pi[\m]$ has finite length
as a smooth $G$-representation.
\end{proposition}

\begin{proof}
Emerton uses geometric Frobenius, so the representation attached
to $\m$ in his convention is $\sigma=\bar\rho^\vee$
\cite[\S2 and Definition~5.2.2]{Em11}. Since
$\sigma\simeq\bar\rho\otimes(\det\bar\rho)^{-1}$, the same local
hypothesis holds for $\sigma_p=\sigma|_{G_{\Q_p}}$.

In the tower with level varying at the primes in $\mathcal R$,
Emerton factors the $\m$-eigenspace as a tensor product of
$\sigma$ and the local representations $\pi(\sigma_v)$ for
$v\in\mathcal R$
\cite[Theorem~6.2.13 and proof of Proposition~6.1.20]{Em11}.
Hochschild--Serre descent at sufficiently deep $p$-level
recovers our fixed tame level by taking $K_{\Sigma_0}$-invariants
\cite[Lemma~5.3.8(2)]{Em11}. Thus
\[
\Pi[\m]\simeq
\sigma\otimes_\F\pi(\sigma_p)\otimes_\F
\left(\bigotimes_{\ell\in\Sigma_0}\pi(\sigma_\ell)\right)^{K_{\Sigma_0}},
\qquad
\sigma_\ell=\sigma|_{G_{\Q_\ell}}.
\]
The tame invariant space is finite-dimensional because each
$\pi(\sigma_\ell)$ is admissible and $K_{\Sigma_0}$ contains a
product of compact open subgroups \cite[Theorem~4.3.1]{Em11}.
Hence $\Pi[\m]$ is a finite direct sum of copies of
$\pi(\sigma_p)$ as a $G$-representation. The latter has finite
length by \cite[Theorem~3.3.2]{Em11}.
\end{proof}

Theorem~\ref{thm:B} now gives noetherianness and $\dim A_\m\ge2$.
For the power series conclusion, we construct a surjection onto
$A_\m$ from the reduction of a fixed-determinant deformation ring.

As in Section~\ref{sec:definite}, let $\Gamma$ be the $p$-Sylow
subgroup of the finite scalar quotient
\[
C_\varepsilon=
\A_f^\times/(\Q^\times K_{Z^p}U_{p,\varepsilon}),
\]
where $K_{Z^p}$ is the subgroup of tame scalars in this $K^p$
and $U_{p,\varepsilon}$ is the kernel of
$z\mapsto\varepsilon(\bar z,\bar z)$ on $\Z_p^\times$.
We use the same fixed-determinant deformation ring
$R^\eta_{\bar\rho,\mathcal R}$, with $\eta\bmod\varpi=\det\bar\rho$.

\begin{corollary}\label{cor:modular-arithmetic}
Assume the local hypothesis of Proposition~\ref{prop:finite-length}.
If $\Gamma=1$ and $H^2(G_{\Q,\mathcal R},\ad^0\bar\rho)=0$, then
\[
A_\m\simeq\F[[X,Y]].
\]
\end{corollary}

\begin{proof}
Let $\mathbb T^{\mathrm{mod}}_\m$ be the inverse limit of the
localized integral tame Hecke images on $H^1(Y_n,\cO)$.
Since $Y_n$ is an open modular curve, $H^2(Y_n,\cO)=0$.
The coefficient sequence therefore gives a Hecke-equivariant
identification
\[
H^1(Y_n,\cO)/\varpi\simeq H^1(Y_n,\F).
\]
Thus $\mathbb T^{\mathrm{mod}}_\m/\varpi$ acts on $\Pi$.
Restricting this action to the spaces $V_{r,\m}$ gives a
continuous map $\mathbb T^{\mathrm{mod}}_\m/\varpi\to A_\m$.

The integral algebra $\mathbb T^{\mathrm{mod}}_\m$ carries a lift
of $\sigma$ with geometric Frobenius polynomial
$X^2-T_\ell X+\ell S_\ell$ \cite[\S5.2]{Em11}.
Taking the dual and base changing to $A_\m$ gives a lift
$\rho:G_{\Q,\mathcal R}\to\GL_2(A_\m)$ of $\bar\rho$ with
\[
\operatorname{tr}\rho(\Frob_\ell)=T_\ell,
\qquad
\det\rho(\Frob_\ell)=\ell S_\ell
\qquad(\ell\notin\mathcal R)
\]
for arithmetic Frobenius.

The scalar action on the $\varepsilon$-spaces factors through
$C_\varepsilon$. Since $\Gamma=1$, each $S_\ell$ has order prime
to $p$ and therefore equals its residue scalar in $\F^\times$
by Hensel's lemma. The Frobenius determinants of $\rho$ and
$\bar\rho$ thus agree, so Chebotarev gives
$\det\rho=\eta\bmod\varpi$. The lift induces a continuous map
\[
R^\eta_{\bar\rho,\mathcal R}/\varpi\longrightarrow A_\m.
\]
Its image is closed by compactness and contains every $T_\ell$
as a Frobenius trace. Since the $S_\ell$ are scalars and
$T_\ell,S_\ell$ topologically generate $A_\m$, the map is
surjective.

As in the proof of Corollary~\ref{cor:arithmetic}, the vanishing
of $H^2$ gives formal smoothness and the global Euler
characteristic formula gives a two-dimensional tangent space.
We obtain a surjection $\F[[X,Y]]\twoheadrightarrow A_\m$.
A nonzero kernel would force $\dim A_\m\le1$, contradicting
the lower bound from Theorem~\ref{thm:B}.
\end{proof}
\section*{Acknowledgements.}
We thank Liang Xiao for introducing us to this topic and for many
helpful discussions. We also thank Jiawei An for his interest and discussions. Part of this work was carried out during a postdoctoral
appointment at the University of Hawai\textquoteleft i at M\=anoa,
whose hospitality and support we gratefully acknowledge.
We thank the contributors to SageMath, which was used for the
computations that motivated this work.

This research was supported by the Beijing Natural Science
Foundation through its International Scientists Project
(Grant No.\ IS26001, 2026--2027).
\section*{AI Declaration}
We used OpenAI's ChatGPT (GPT-6) to assist
with editing and the organization of this paper, and
to help check calculations and review the arguments.
We reviewed all AI-assisted revisions and take full
responsibility for the contents of this paper.
\renewcommand{\bibliofont}{\small}

\end{document}